\documentclass[12pt]{article}

\usepackage{fullpage}
\usepackage{enumitem}
\usepackage{lineno}
\usepackage[b]{esvect}
\usepackage{tikz}
\usepackage{caption}
\usepackage{subcaption}
\usepackage[font=footnotesize,labelfont=bf]{caption}
\usetikzlibrary{shapes.multipart,shapes.geometric,topaths,calc}
\usepackage{epsfig}
\usepackage{amsfonts}
\usepackage{amssymb}
\usepackage{amstext}
\usepackage{amsmath}
\usepackage{xspace}
\usepackage{theorem}
\usepackage{xcolor}
\usepackage{color}
\usepackage{graphicx}
\usepackage{ifthen}
\usepackage{hyperref}
\usetikzlibrary{shapes,arrows}

\newcommand{\junk}[1]{}

\makeatletter
\newenvironment{proof}{{\bf Proof:  }}{\hfill\rule{2mm}{2mm}}

\newtheorem{theorem}{Theorem}
\newtheorem{lemma}[theorem]{Lemma}
\newtheorem{proposition}[theorem]{Proposition}

\newcommand{\sq}{~\square~}
\newcommand{\F}{{\cal F}}

\title{Decompositions of even hypercubes into cycles whose length is a power of two}

\author{
Samuel Gibson\thanks{Carnegie Mellon University, \texttt{sgibson@andrew.cmu.edu}} 
\and David Offner\thanks{Carnegie Mellon University, \texttt{doffner@andrew.cmu.edu}}
}

\date{\today}

\begin{document}

\maketitle

\begin{abstract}
If $n$ is even, the $n$-dimensional hypercube can be decomposed into edge-disjoint cycles of length $2^i$ for every value of $i$ from $2$ to $n$.
\end{abstract}

\section{Introduction}

A decomposition of a graph $G$ is a set of subgraphs that partition the edges of $G$. For $n \ge 1$, the $n$-dimensional hypercube graph is denoted by $Q_n$ and is defined as the graph with $V(Q_n) = \{0,1\}^n$, and 
\[E(Q_n) = \{ vw : v,w \in V(Q_n), \text{$v$ and $w$ differ in one coordinate}\}.\]
We are interested in decomposing the hypercube into cycles each of whose length is the same power of two.


Decomposition problems on the hypercube go back to Ringel~\cite{ringel}, who proved that $Q_n$ has a decomposition into Hamiltonian cycles when $n$ is a power of two and $n \ge 2$. Such a decomposition of a graph into Hamiltonian cycles is called a Hamiltonian decomposition, and Ringel asked whether $Q_n$ has a Hamiltonian decomposition for all even $n$. This was implicitly resolved in the affirmative by Aubert and Schneider~\cite{aubert}. See Alspach, Bermond, and Sotteau~\cite{alspach} for an explicit statment. 

Beyond Hamiltonian decompositions, Stout~\cite{stout}, Horak, Siran, and Wallis~\cite{HSW}, Mollard and Ramras~\cite{MR}, and Wagner and Wild \cite{WW12} showed that $Q_n$ can be decomposed into certain trees.
Anick and Ramras~\cite{AR}, and independently Erde~\cite{erde} proved that if $n$ is odd, then $Q_n$ can be decomposed into any path whose length divides the number of edges in $Q_n$ and is at most $n$.
Fink~\cite{fink}, Horak, Siran, and Wallis~\cite{HSW}, Mollard and Ramras~\cite{MR}, Tapadia, Waphare, and Borse~\cite{TWB}, and Axenovich, Offner, and Tompkins~\cite{AOT21} proved that when $n$ is even, $Q_n$ can be decomposed into cycles of certain lengths, among other results, but a complete characterization of the graphs that can decompose $Q_n$ remains an open question, even for paths, trees, and cycles.

Our goal is to find decompositions of hypercube graphs into cycles whose length is a power of two.  The results of Ringel~\cite{ringel}, Aubert and Schneider~\cite{aubert}, and Axenovich, Offner, and Tompkins~\cite{AOT21} show that $Q_n$ can be decomposed into relatively long cycles whose length is a power of two, and  Tapadia, Waphare, and Borse ~\cite{TWB} show that $Q_n$ can be decomposed into relatively short cycles whose length is a power of two, but there are no published  results for cycles whose length is in between.  In this paper we give an elementary proof of the following complete characterization.

\begin{theorem}\label{alldecomp}
If $n$ is even, $n \ge 2$, and $2 \le i \le n$, then the $n$-dimensional hypercube can be decomposed into cycles of length $2^i$.
\end{theorem}

Theorem~\ref{alldecomp} is a corollary of our main result, the slightly more general Theorem~\ref{partitionable}.  In Section~\ref{defs}, we introduce the definitions and notation necessary for our proofs.   In Section~\ref{secdec} we prove the main result.  In Proposition~\ref{halfham} we prove that if $n$ is even, then $Q_n$ can be decomposed into cycles of length $2^{n-1}$, which, along with Theorem~\ref{ham}, provides the base cases for the inductive proof of Theorem~\ref{partitionable}.  For the inductive step, in Lemma~\ref{8ell} and Proposition~\ref{cart} we use the fact that $Q_{n+2}$ is the Cartesian product of $Q_n$ with $Q_2$ to show that if $Q_n$ can be decomposed into cycles of length $2^i$, then $Q_{n+2}$ can as well.

\section{Definitions and notation}\label{defs}

The Cartesian product of two graphs $G$ and $H$ is denoted as $G \sq H$ and is defined to be the graph where $V(G \sq H) = \{ (u,v): u \in V(G), v \in V(H)\}$ and $E(G \sq H) = \{(u, v)(u', v') : u = u', vv' \in E(H) \text{ or } v = v', uu' \in E(G)\}$.  Define the $n$-fold Cartesian product of a graph $G$ with itself to be the graph $G^n = G \sq G \sq \cdots \sq G$ resulting from taking the Cartesian product of $G$ with itself $n$ times.  In other words, $V(G^n) = V(G)^n$ and the edges of $G^n$ are pairs $\{(v_1, \ldots, v_n),(w_1, \ldots, w_n)\}$ where  $(v_1, \ldots, v_n) ,(w_1, \ldots, w_n) \in V(G^n)$, and $(v_1, \ldots, v_n)$ and $(w_1, \ldots, w_n)$ differ in only coordinate $i$, where $v_iw_i \in E(G)$.

Let $C_k$ represent the cycle of length $k$. Let $K_n$ represent the complete graph on $n$ vertices. We call the Cartesian product of two cycles a torus.

For $n \ge 1$, the $n$-dimensional hypercube graph can be represented as a Cartesian product.
Since we can represent $K_2$ as the graph with $V(K_2) = \{0,1\}$, and $E(K_2) = \{\{0,1\}\}$, $Q_1 = K_2$, and for $n \ge 2$, $Q_n = K_2^n$ is the $n$-fold Cartesian product of $K_2$. Since $Q_2 =  K_2 \sq K_2=C_4$, for $n \ge 1$, $Q_{2n} = C_4^n$ is the $n$-fold Cartesian product of $C_4$.

\begin{figure}
    \begin{tikzpicture}[scale = .8]

    \draw[ultra thick] (1.5, 2) -- (2,2) -- (2,3) -- (3, 3) -- (3, 2) -- (4,2) -- (4,3) -- (5,3) -- (5,2) -- (5.5, 2);
    
    \draw[ultra thick] (1.5, 3) -- (2,3) -- (2,4) -- (3, 4) -- (3, 3) -- (4,3) -- (4,4) -- (5,4) -- (5,3) -- (5.5, 3);
    
    \draw[ultra thick] (1.5, 4) -- (2,4) -- (2,5) -- (3, 5) -- (3, 4) -- (4,4) -- (4,5) -- (5,5) -- (5,4) -- (5.5, 4);
    
    \draw[ultra thick] (1.5, 5) -- (2,5) -- (2, 5.5);
    
    \draw[ultra thick] (3, 5.5) -- (3,5) -- (4, 5) -- (4, 5.5);
    
    \draw[ultra thick] (5, 5.5) -- (5,5) -- (5.5, 5);
    
    \draw[ultra thick] (2, 1.5) -- (2,2) -- (3,2) -- (3, 1.5);
    
    \draw[ultra thick] (4, 1.5) -- (4,2) -- (5,2) -- (5, 1.5);

    \foreach \x in {2,3,4, 5}{
    \foreach \y in {2,3,4, 5}{
    \node[style=circle, draw, fill=black,scale=.4] at (\x,\y){};}};
    \node at (1,5){00};
    \node at (1,4){01};
    \node at (1,3){11};
    \node at (1,2){10};
    \node at (2,6){00};
    \node at (3,6){01};
    \node at (4,6){11};
    \node at (5,6){10};
    
 \node[style=circle, draw, fill=white,scale=.6] at (4,2){};

\end{tikzpicture}
\hfill
    \begin{tikzpicture}[scale = .8]

    \draw[ultra thick, red] (1.5, 2) -- (2,2) -- (2,3) -- (3, 3) -- (3, 2) -- (4,2) -- (4,3) -- (5,3) -- (5,2) -- (5.5, 2);
    
    \draw[ultra thick, dashed, blue] (1.5, 3) -- (2,3) -- (2,4) -- (3, 4) -- (3, 3) -- (4,3) -- (4,4) -- (5,4) -- (5,3) -- (5.5, 3);
    
    \draw[ultra thick, red] (1.5, 4) -- (2,4) -- (2,5) -- (3, 5) -- (3, 4) -- (4,4) -- (4,5) -- (5,5) -- (5,4) -- (5.5, 4);
    
    \draw[ultra thick, dashed, blue] (1.5, 5) -- (2,5) -- (2, 5.5);
    
    \draw[ultra thick, dashed, blue] (3, 5.5) -- (3,5) -- (4, 5) -- (4, 5.5);
    
    \draw[ultra thick, dashed, blue] (5, 5.5) -- (5,5) -- (5.5, 5);
    
    \draw[ultra thick, dashed, blue] (2, 1.5) -- (2,2) -- (3,2) -- (3, 1.5);
    
    \draw[ultra thick, dashed, blue] (4, 1.5) -- (4,2) -- (5,2) -- (5, 1.5);

    \foreach \x in {2,3,4, 5}{
    \foreach \y in {2,3,4, 5}{
    \node[style=circle, draw, fill=black,scale=.4] at (\x,\y){};}};
    \node at (1,5){00};
    \node at (1,4){01};
    \node at (1,3){11};
    \node at (1,2){10};
    \node at (2,6){00};
    \node at (3,6){01};
    \node at (4,6){11};
    \node at (5,6){10};
\end{tikzpicture}
\hfill
\begin{tikzpicture}[scale = .8]

    \draw[ultra thick, red] (1.5, 2) -- (2,2) -- (2,3) -- (3, 3) -- (3, 2) -- (4,2) -- (4,3) -- (5,3) -- (5,4) -- (5.5,4);
    
    \draw[ultra thick, red] (1.5, 4) -- (2, 4) -- (2,5) -- (3,5) -- (3,4) -- (4,4) -- (4,5) -- (5,5) -- (5, 5.5);
    
    \draw[ultra thick, dashed, blue] (1.5, 3) -- (2,3) -- (2,4) -- (3, 4) -- (3, 3) -- (4,3) -- (4,4) -- (5,4) -- (5,5) -- (5.5, 5); 
    
    \draw[ultra thick, dashed, blue] (5,3) -- (5.5, 3);

    \draw[ultra thick, dashed, blue] (1.5, 5) -- (2,5) -- (2, 5.5);
    
    \draw[ultra thick, dashed, blue] (3, 5.5) -- (3,5) -- (4, 5) -- (4, 5.5);
    
    
    \draw[ultra thick, dashed, blue] (2, 1.5) -- (2,2) -- (3,2) -- (3, 1.5);
    
    \draw[ultra thick, dashed, blue] (4, 1.5) -- (4,2) -- (5,2) -- (5, 3);
    
    \draw[ultra thick, red](5, 1.5) -- (5,2)  -- (5.5, 2);

    \foreach \x in {2,3,4, 5}{
    \foreach \y in {2,3,4, 5}{
    \node[style={circle, draw, fill=black},scale=.4] at (\x,\y){};}};
    
    \node at (1,5){00};
    \node at (1,4){01};
    \node at (1,3){11};
    \node at (1,2){10};
    \node at (2,6){00};
    \node at (3,6){01};
    \node at (4,6){11};
    \node at (5,6){10};
    
\end{tikzpicture}
\caption{
Left: A representation of $Q_4$ as the Cartesian product $C_4 \sq C_4$.  With the vertices of $C_4$ labeled $00, 01, 11, 10$, the vertex labels for  $Q_4$ are obtained by concatenating the column label of the vertex with the row label of the vertex. For example the large white vertex has label 1110.
Center: A partitionable decomposition of $Q_4$ into cycles of length 8. One partition set $\F_1$ contains the two solid red cycles, while the other partition set $\F_2$ contains the two dashed blue cycles.
Right: A partitionable decomposition of $Q_4$ into cycles of length 16. One partition set $\F_1$ contains the solid red cycle, while the other partition set $\F_2$ contains the dashed blue cycle.}
\label{Q_4}
\end{figure}
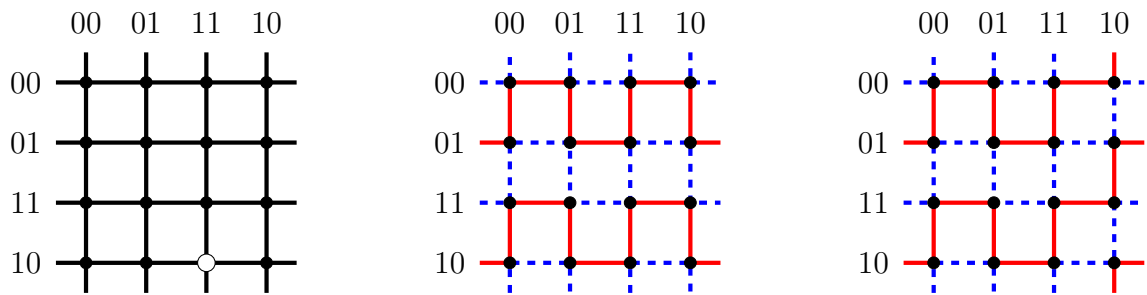

See Figure~\ref{Q_4} (left) for a representation of $Q_4$ as the Cartesian product $C_4 \sq C_4$.  With the vertices of each $C_4$ labeled $00, 01, 11, 10$, the vertex labels for the vertices of $Q_4$ are obtained by concatenating the column label of the vertex with the row label of the vertex.  For a Cartesian product $G \sq H$, we follow the convention in Figure~\ref{Q_4} and refer to the edges $\{(u, v)(u', v') : u = u', vv' \in E(H)\}$ as vertical edges, and the edges $\{(u, v)(u', v') : v = v', uu' \in E(G)\}$ as horizontal edges. In other words, an edge whose vertices differ in the first coordinate is referred to as a horizontal edge, and an edge whose vertices differ in the second coordinate is a vertical edge.  

Given two graphs $G$ and $H$, define the union $G \cup H$ to be the graph where $V(G \cup H) = V(G) \cup V(H)$ and $E(G \cup H)= E(G) \cup E(H)$. We use the notation $G \sqcup H$ to denote an edge-disjoint union of graphs, that is a union of graphs $G$ and $H$ where $E(G) \cap E(H) = \emptyset$.

A decomposition of a graph $G$ is a set of subgraphs that partition the edges in $G$. That is, a set of graphs $\{G_1, G_2, ..., G_m \}$ is a decomposition of $G$ if and only if $E(G) = \bigsqcup\limits_{1 \leq i \leq m} E(G_i)$. Note $ E(G_i) \cap E(G_j) = \emptyset$ if $i \neq j$. We will write $G = G_1 \sqcup \cdots\sqcup G_m$ to represent such a decomposition. 

Call a decomposition $G = G_1 \sqcup \cdots \sqcup G_m$ partitionable if the set $\{G_1, \ldots, G_m\}$ can be partitioned into subsets $\F_1, \ldots, \F_k$ such that for $1 \le i \le k$ the sets $\{V(H) : H \in \F_i\}$ partition the vertex set of $G$. Call the sets $\F_1, \ldots, \F_k$ the partition sets of the decomposition.  For example, 
Figure~\ref{Q_4} (center) shows a partitionable decomposition of $Q_4$ into cycles of length 8.  Figure~\ref{Q_4} (right) shows a partitionable decomposition of $Q_4$ into cycles of length 16.

\section{Decompositions}\label{secdec}

Propositions~\ref{q4} and \ref{q6} describe decompositions of $Q_4$ and $Q_6$ into cycles.  These decompositions serve as base cases for Proposition~\ref{halfham}. 

\begin{proposition}\label{q4}
$Q_4$ has a partitionable decomposition into four cycles of length $8$.
\end{proposition}

\begin{proof}
The decomposition is shown in Figure~\ref{Q_4} (center). The two red cycles are $R^0$ = (0000, 0100, 0101, 1101, 1100, 1000, 1001, 0001, 0000) and  $R^1 =$ (0011, 0111, 0110, 1110, 1111, 1011, 1010, 0010, 0011) and the two dashed blue cycles are $B^0 =$ (0000, 0010, 0110, 0100, 1100, 1110, 1010, 1000, 0000) and $B^1 =$ (0011, 0001, 0101, 0111, 1111, 1101, 1001, 1011, 0011). The two partition sets are $\{R^0, R^1\}$ and $\{B^0, B^1\}$.
\end{proof}

\begin{proposition}\label{q6}
$Q_6$ has a partitionable decomposition into six cycles of length $32$.
\end{proposition}

\begin{figure}
\begin{tikzpicture}[scale = 1.3]

    \draw[gray] (0.5, 12.5) -- (5.9, 12.5) -- (5.9, 7.2) -- (0.5, 7.2) -- (0.5,12.5);
    
    \draw[gray] (7.1, 12.5) -- (12.5, 12.5) -- (12.5, 7.2) -- (7.1, 7.2) -- (7.1,12.5);
    
    \draw[gray] (7.1, 5.75) -- (12.5, 5.75) -- (12.5, .45) -- (7.1, .45) -- (7.1,5.75);
    
    \draw[gray] (.5, 5.75) -- (5.9, 5.75) -- (5.9, .45) -- (.5, .45) -- (.5,5.75);

    \node[font = \bfseries \LARGE] at (1,12){00};
    \node[font = \LARGE \bfseries] at (12,12){01};
    \node[font = \LARGE \bfseries] at (12,1){11};
    \node[font = \LARGE \bfseries] at (1,1){10};

    
    \node at (2,1){00};
    \node at (3,1){01};
    \node at (4,1){11};
    \node at (5,1){10};
    
    \node at (1,5){00};
    \node at (1,4){01};
    \node at (1,3){11};
    \node at (1,2){10};

    \draw (5, 5) -- (5, 5.5);
    
    \draw[ultra thick, red, dashed] (1.5, 2) -- (2,2) -- (2,3) -- (3, 3) -- (3, 2);
    \draw[ultra thick, red, dashed] (4,2) -- (4,3) -- (5,3) -- (5,2) -- (5.5, 2);
    
    \draw[ultra thick, red, dashed] (3,2) to[out=-78,in=-72] (3,8);
    
    \draw[ultra thick, blue, dotted] (3,4) to[out=250,in=250] (3,10);
    
    \draw[ultra thick, gray] (1.5, 3) -- (2,3) -- (2,4) -- (3, 4) -- (3, 3) -- (4,3) -- (4,4) -- (5,4) -- (5,3) -- (5.5, 3);
    
    \draw[ultra thick, brown] (3,4) -- (4,4);
     \draw [ultra thick, yellow](3, 2) -- (4,2);
    
    \draw[ultra thick, blue, dotted] (1.5, 4) -- (2,4) -- (2,5) -- (3, 5) -- (3, 4);
    \draw[ultra thick, blue, dotted] (4,4) -- (4,5) -- (5,5) -- (5,4) -- (5.5, 4);
    
    \draw[ultra thick, gray] (1.5, 5) -- (2,5) -- (2, 5.5);
    
    \draw[ultra thick, gray] (3, 5.5) -- (3,5) -- (4, 5) -- (4, 5.5);
    
    \draw[ultra thick, gray] (5, 5.5) -- (5,5) -- (5.5, 5);
    
    \draw[ultra thick, gray] (2, 1.5) -- (2,2) -- (3,2) -- (3, 1.5);
    
    \draw[ultra thick, gray] (4, 1.5) -- (4,2) -- (5,2) -- (5, 1.5);
    
    \foreach \x in {2,3,4, 5}{
    \foreach \y in {2,3,4, 5}{
    \node[style=circle, draw, fill=black,scale=.4] at (\x,\y){};}};


    \node at (8,1){00};
    \node at (9,1){01};
    \node at (10,1){11};
    \node at (11,1){10};

    \node at (12,5){00};
    \node at (12,4){01};
    \node at (12,3){11};
    \node at (12,2){10};
    

    \draw (11, 5) -- (11, 5.5);
    
    \draw[ultra thick, red, dashed] (7.5, 2) -- (8,2) -- (8,3) -- (9, 3) -- (9, 2);
    \draw[ultra thick, red, dashed] (10,2) -- (10,3) -- (11,3) -- (11,2) -- (11.5, 2);
    
    \draw[ultra thick, gray] (7.5, 3) -- (8,3) -- (8,4) -- (9, 4) -- (9, 3) -- (10,3) -- (10,4) -- (11,4) -- (11,3) -- (11.5, 3);
    
    \draw[ultra thick, blue, dotted] (7.5, 4) -- (8,4) -- (8,5) -- (9, 5) -- (9, 4);
    
    \draw[ultra thick, blue, dotted] (10,4) to[out=195,in=210] (4,4);
    
    \draw[ultra thick, red, dashed] (10,2) to[out=195,in=210] (4,2);
    
    \draw[ultra thick, blue, dotted](10,4) -- (10,5) -- (11,5) -- (11,4) -- (11.5, 4);
    
    \draw[ultra thick, gray] (7.5, 5) -- (8,5) -- (8, 5.5);
    
    \draw[ultra thick, gray] (9, 5.5) -- (9,5) -- (10, 5) -- (10, 5.5);
    
    \draw[ultra thick, gray] (11, 5.5) -- (11,5) -- (11.5, 5);
    
    \draw[ultra thick, gray] (8, 1.5) -- (8,2) -- (9,2) -- (9, 1.5);
    
    \draw [ultra thick, yellow](9, 2) -- (10,2);
    
    \draw[ultra thick, brown] (9,4) -- (10,4);
    
    \draw[ultra thick, gray] (10, 1.5) -- (10,2) -- (11,2) -- (11, 1.5);

    \foreach \x in {8,9,10, 11}{
    \foreach \y in {2,3,4, 5}{
    \node[style=circle, draw, fill=black,scale=.4] at (\x,\y){};}};

    
    \node at (8,12){00};
    \node at (9,12){01};
    \node at (10,12){11};
    \node at (11,12){10};

    \node at (12,11){00};
    \node at (12,10){01};
    \node at (12,9){11};
    \node at (12,8){10};
    
    \draw[ultra thick, yellow] (9, 8) -- (10,8);
     \draw [ultra thick, brown](9, 10) -- (10,10);
    
    \draw (11, 11) -- (11, 11.5);
    
    \draw[ultra thick, red, dashed] (7.5, 8) -- (8,8) -- (8,9) -- (9, 9) -- (9, 8);
    
    \draw[ultra thick, red,  dashed] (9,8) to[out=110,in=110] (9,2);
    
    \draw[ultra thick, blue,  dotted] (9,10) to[out=70,in=70] (9,4);
    
    \draw[ultra thick, red, dashed](10,8) -- (10,9) -- (11,9) -- (11,8) -- (11.5, 8);
    
    \draw[ultra thick, gray] (7.5, 9) -- (8,9) -- (8,10) -- (9, 10) -- (9, 9) -- (10,9) -- (10,10) -- (11,10) -- (11,9) -- (11.5, 9);
    
    \draw[ultra thick, blue, dotted] (7.5, 10) -- (8,10) -- (8,11) -- (9, 11) -- (9, 10);
    \draw[ultra thick, blue, dotted](10,10) -- (10,11) -- (11,11) -- (11,10) -- (11.5, 10);
    
    \draw[ultra thick, gray] (7.5, 11) -- (8,11) -- (8, 11.5);
    
    \draw[ultra thick, gray] (9, 11.5) -- (9,11) -- (10, 11) -- (10, 11.5);
    
    \draw[ultra thick, gray] (11, 11.5) -- (11,11) -- (11.5, 11);
    
    \draw[ultra thick, gray] (8, 7.5) -- (8,8) -- (9,8) -- (9, 7.5);
    
    \draw[ultra thick, gray] (10, 7.5) -- (10,8) -- (11,8) -- (11, 7.5);

    \foreach \x in {8,9,10, 11}{
    \foreach \y in {8,9,10, 11}{
    \node[style=circle, draw, fill=black,scale=.4] at (\x,\y){};}};

    
    \node at (1,11){00};
    \node at (1,10){01};
    \node at (1,9){11};
    \node at (1,8){10};
    
    \node at (2,12){00};
    \node at (3,12){01};
    \node at (4,12){11};
    \node at (5,12){10};

    \draw[ultra thick, brown] (3,10) -- (4,10);
    \draw[ultra thick, yellow] (3,8) -- (4,8);
    
    \draw (5, 11) -- (5, 11.5);
    
    \draw[ultra thick, red, dashed] (1.5, 8) -- (2,8) -- (2,9) -- (3, 9) -- (3, 8); 
    
    \draw[ultra thick, red,  dashed] (4,8) to[out=35,in=15] (10,8);
    
    \draw[ultra thick, blue,  dotted] (4,10) to[out=35,in=15] (10,10);
    
    \draw[ultra thick, red, dashed] (4,8) -- (4,9) -- (5,9) -- (5,8) -- (5.5, 8);
    
    \draw[ultra thick, gray] (1.5, 9) -- (2,9) -- (2,10) -- (3, 10) -- (3, 9) -- (4,9) -- (4,10) -- (5,10) -- (5,9) -- (5.5, 9);
    
    \draw[ultra thick, blue, dotted] (1.5, 10) -- (2,10) -- (2,11) -- (3, 11) -- (3, 10);
    \draw[ultra thick, blue, dotted](4,10) -- (4,11) -- (5,11) -- (5,10) -- (5.5, 10);
    
    \draw[ultra thick, gray] (1.5, 11) -- (2,11) -- (2, 11.5);
    
    \draw[ultra thick, gray] (3, 11.5) -- (3,11) -- (4, 11) -- (4, 11.5);
    
    \draw[ultra thick, gray] (5, 11.5) -- (5,11) -- (5.5, 11);
    
    \draw[ultra thick, gray] (2, 7.5) -- (2,8) -- (3,8) -- (3, 7.5);
    
    \draw[ultra thick, gray] (4, 7.5) -- (4,8) -- (5,8) -- (5, 7.5);

    \foreach \x in {2,3,4, 5}{
    \foreach \y in {8,9,10, 11}{
    \node[style=circle, draw, fill=black,scale=.4] at (\x,\y){};}};
    
    \node[style=circle, draw, fill=white,scale=.6] at (2,3){};
    
\end{tikzpicture}
\caption{A representation of $Q_6 = Q_4 \sq C_4$. The  $Q_4$ whose vertex labels end with 00, 01, 11, and 10 are pictured in the upper left, upper right, lower right, and lower left, respectively. The label for each vertex in $Q_6$ is obtained by concatenating the column label and then the row label of the vertex, followed by these last two coordinates. For example, the address of the large white vertex in the lower left is 001110.  All horizontal edges (edges that differ in the first four coordinates) are shown, but the eight curved edges are the only vertical edges shown, and the rest are omitted for clarity.  The dashed red and dotted blue cycles are $C^1$ and $C^2$ in the decomposition of $Q_6$ described in the proof of Proposition~\ref{q6}.  The remaining gray, brown, yellow, and omitted vertical edges are partitioned between the graphs $G_B$ and $G_Y$ in Figure~\ref{Q_6 C_32 Decomp} (top and middle).  The brown edges are $\{010100, 110100\}$, $\{010110, 110110\}$, $\{010111, 110111\}$, and $\{010101, 110101\}$. The yellow edges are $\{011000, 111000\}$, $\{011010, 111010\}$, $\{011011, 111011\}$, and $\{011001, 111001\}$.
}
\label{Q_6C_32}
\end{figure}
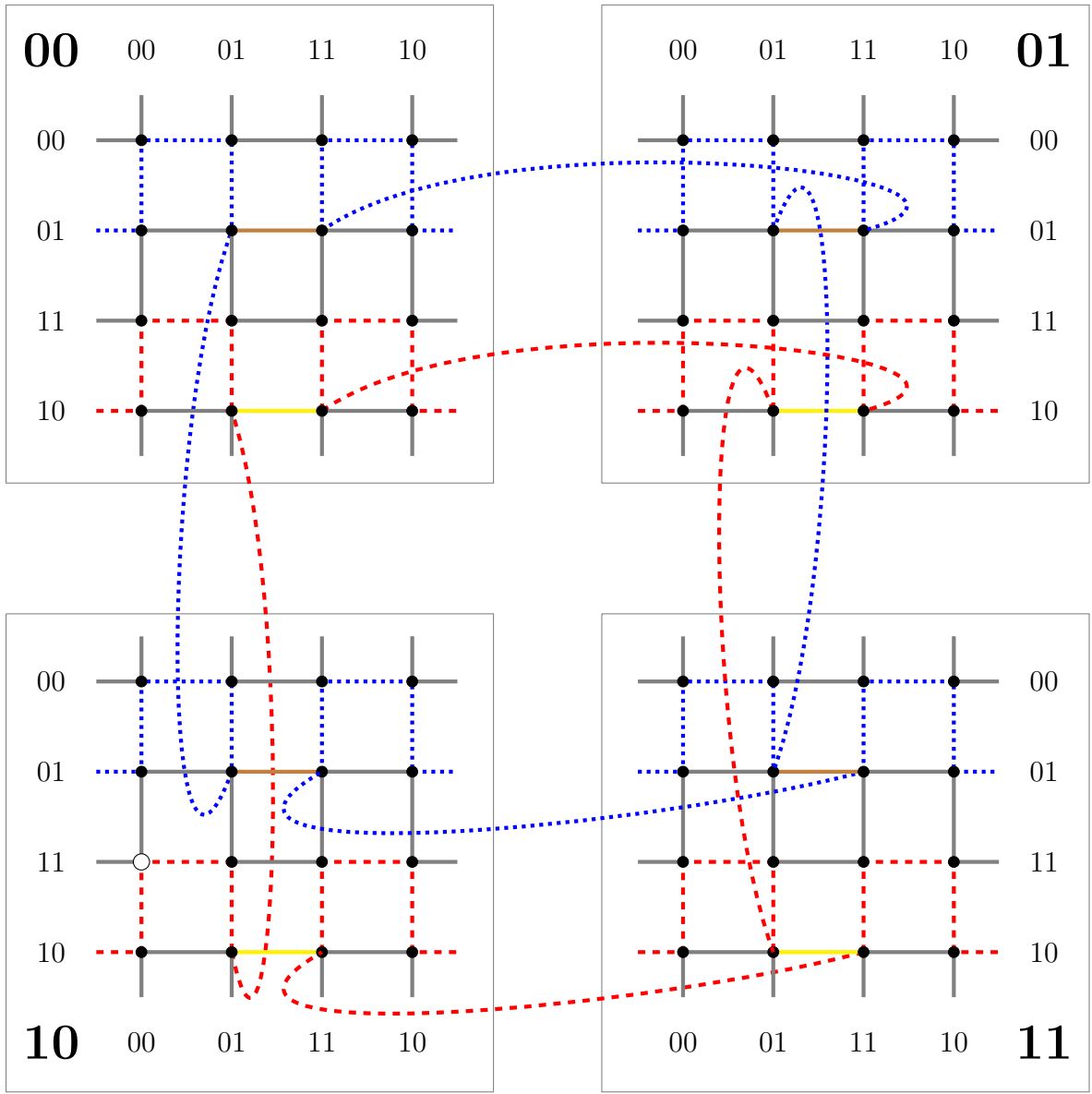

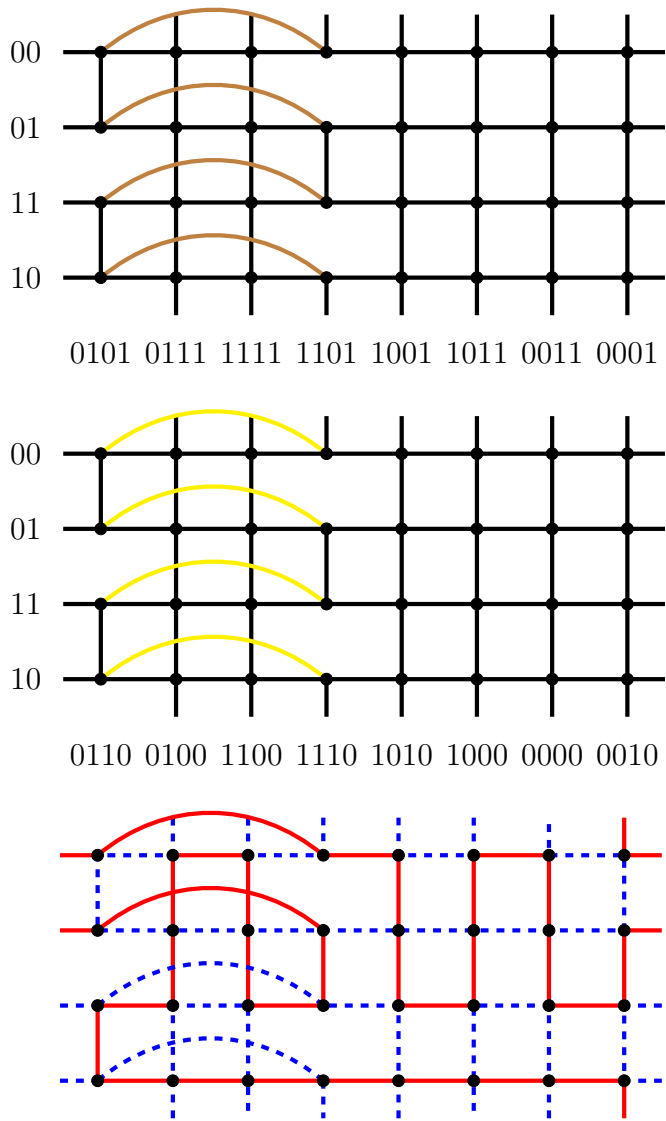
\begin{figure}
\begin{center}
    \begin{tikzpicture}[scale = 1]

\node at (0,-1){0101};
\node at (1,-1){0111};
\node at (2,-1){1111};
\node at (3,-1){1101};
\node at (4,-1){1001};
\node at (5,-1){1011};
\node at (6,-1){0011};
\node at (7,-1){0001};

\node at (-1,3){00};
\node at (-1,2){01};
\node at (-1,1){11};
\node at (-1,0){10};

\draw[ultra thick] (-.5,3) -- (0,3);
\draw[ultra thick] (-.5,2) -- (0,2);
\draw[ultra thick] (7.5,3) -- (7,3) -- (7,3.5);
\draw[ultra thick] (7, -.5) -- (7,0) -- (0,0) -- (0,1) -- (1,1) -- (1,3) -- (2,3) -- (2,1) -- (3,1) -- (3,2);
\draw[ultra thick] (3,3) -- (4,3) -- (4,1) -- (5,1) -- (5,3) -- (6,3) -- (6,1) -- (7,1) -- (7,2) -- (7.5,2);

\draw[ultra thick] (1,3.5) -- (1,3) -- (0,3) -- (0,2) -- (7,2) -- (7,3) -- (6,3) -- (6,3.5);

\draw[ultra thick] (1,-.5) -- (1,1) -- (2,1) -- (2,-.5);
\draw[ultra thick] (6,-.5) -- (6,1) -- (5,1) -- (5,-.5);
\draw[ultra thick] (5,3.5) -- (5,3) -- (4,3) -- (4,3.5);
\draw[ultra thick] (4,-.5) -- (4,1) -- (3,1);
\draw[ultra thick] (2,3.5) -- (2,3) -- (3,3) -- (3, 3.5);
\draw[ultra thick] (-0.5,1) -- (0,1);
\draw[ultra thick] (-0.5,0) -- (0,0);
\draw[ultra thick] (7.5,1) -- (7,1) -- (7,0) -- (7.5,0);
\draw[ultra thick] (3,0) -- (3,-.5);
\draw[ultra thick,  brown] (3,2) to[out=140,in=40] (0,2);

\draw[ultra thick, brown] (3,3) to[out=140,in=40] (0,3);

\draw[ultra thick, brown] (3,1) to[out=140,in=40] (0,1);

\draw[ultra thick, brown] (3,0) to[out=140,in=40] (0,0);

\foreach \x in {0, ..., 7}{
    \foreach \y in {0,...,3}{
    \node[style={circle, draw, fill=black},scale=.4] at (\x,\y){};}};

\node at (-1.1,3){};

\end{tikzpicture}

\vspace{.1in}

\begin{tikzpicture}[scale = 1]

\node at (0,-1){0110};
\node at (1,-1){0100};
\node at (2,-1){1100};
\node at (3,-1){1110};
\node at (4,-1){1010};
\node at (5,-1){1000};
\node at (6,-1){0000};
\node at (7,-1){0010};

\node at (-1,3){00};
\node at (-1,2){01};
\node at (-1,1){11};
\node at (-1,0){10};

\draw[ultra thick] (-.5,3) -- (0,3);
\draw[ultra thick] (-.5,2) -- (0,2);
\draw[ultra thick] (7.5,3) -- (7,3) -- (7,3.5);
\draw[ultra thick] (7, -.5) -- (7,0) -- (0,0) -- (0,1) -- (1,1) -- (1,3) -- (2,3) -- (2,1) -- (3,1) -- (3,2);
\draw[ultra thick] (3,3) -- (4,3) -- (4,1) -- (5,1) -- (5,3) -- (6,3) -- (6,1) -- (7,1) -- (7,2) -- (7.5,2);

\draw[ultra thick] (1,3.5) -- (1,3) -- (0,3) -- (0,2) -- (7,2) -- (7,3) -- (6,3) -- (6,3.5);

\draw[ultra thick] (1,-.5) -- (1,1) -- (2,1) -- (2,-.5);
\draw[ultra thick] (6,-.5) -- (6,1) -- (5,1) -- (5,-.5);
\draw[ultra thick] (5,3.5) -- (5,3) -- (4,3) -- (4,3.5);
\draw[ultra thick] (4,-.5) -- (4,1) -- (3,1);
\draw[ultra thick] (2,3.5) -- (2,3) -- (3,3) -- (3, 3.5);
\draw[ultra thick] (-0.5,1) -- (0,1);
\draw[ultra thick] (-0.5,0) -- (0,0);
\draw[ultra thick] (7.5,1) -- (7,1) -- (7,0) -- (7.5,0);
\draw[ultra thick] (3,0) -- (3,-.5);
\draw[ultra thick, yellow] (3,2) to[out=140,in=40] (0,2);

\draw[ultra thick, yellow] (3,3) to[out=140,in=40] (0,3);

\draw[ultra thick, yellow] (3,1) to[out=140,in=40] (0,1);

\draw[ultra thick, yellow] (3,0) to[out=140,in=40] (0,0);

\foreach \x in {0, ..., 7}{
    \foreach \y in {0,...,3}{
    \node[style={circle, draw, fill=black},scale=.4] at (\x,\y){};}};
    
\node at (-1.1,3){};

\end{tikzpicture}

\vspace{.1in}

\begin{tikzpicture}[scale = 1]
\draw[ultra thick, red] (-.5,3) -- (0,3);
\draw[ultra thick, red] (-.5,2) -- (0,2);
\draw[ultra thick, red] (7.5,3) -- (7,3) -- (7,3.5);
\draw[ultra thick, red] (7, -.5) -- (7,0) -- (0,0) -- (0,1) -- (1,1) -- (1,3) -- (2,3) -- (2,1) -- (3,1) -- (3,2);
\draw[ultra thick, red] (3,3) -- (4,3) -- (4,1) -- (5,1) -- (5,3) -- (6,3) -- (6,1) -- (7,1) -- (7,2) -- (7.5,2);

\draw[ultra thick, blue, dashed] (1,3.5) -- (1,3) -- (0,3) -- (0,2) -- (7,2) -- (7,3) -- (6,3) -- (6,3.5);

\draw[ultra thick, blue, dashed] (1,-.5) -- (1,1) -- (2,1) -- (2,-.5);
\draw[ultra thick, blue, dashed] (6,-.5) -- (6,1) -- (5,1) -- (5,-.5);
\draw[ultra thick, blue, dashed] (5,3.5) -- (5,3) -- (4,3) -- (4,3.5);
\draw[ultra thick, blue, dashed] (4,-.5) -- (4,1) -- (3,1);
\draw[ultra thick, blue, dashed] (2,3.5) -- (2,3) -- (3,3) -- (3, 3.5);
\draw[ultra thick, blue, dashed] (-0.5,1) -- (0,1);
\draw[ultra thick, blue, dashed] (-0.5,0) -- (0,0);
\draw[ultra thick, blue, dashed] (7.5,1) -- (7,1) -- (7,0) -- (7.5,0);
\draw[ultra thick, blue, dashed] (3,0) -- (3,-.5);
\draw[ultra thick, red] (3,2) to[out=140,in=40] (0,2);

\draw[ultra thick, red] (3,3) to[out=140,in=40] (0,3);

\draw[ultra thick, blue, dashed] (3,1) to[out=140,in=40] (0,1);

\draw[ultra thick, blue, dashed] (3,0) to[out=140,in=40] (0,0);

\foreach \x in {0, ..., 7}{
    \foreach \y in {0,...,3}{
    \node[style={circle, draw, fill=black},scale=.4] at (\x,\y){};}};
    
    \node at (-1.1,3){};

\end{tikzpicture}
\end{center}
\captionof{figure}{Top: The subgraph $G_B$ of $Q_6$. The curved brown edges are the same as the brown edges in Figure~\ref{Q_6C_32}. The black edges are either gray edges or omitted vertical edges in Figure~\ref{Q_6C_32}. The four vertical ``missing'' edges in this figure are the dashed blue curved edges in Figure~\ref{Q_6C_32}. Middle: The subgraph $G_Y$ of $Q_6$. The curved yellow edges are the same as the yellow edges in Figure~\ref{Q_6C_32}. The black edges are either gray edges or omitted vertical edges in Figure~\ref{Q_6C_32}. The four vertical ``missing'' edges in this figure are the dashed curved red edges in Figure~\ref{Q_6C_32}. Bottom: A Hamiltonian decomposition for $G_B$ or $G_Y$.}
\label{Q_6 C_32 Decomp}
\end{figure}

\begin{proof}
A partitionable decomposition of $Q_6$ into six cycles of length $32$ is illustrated in Figures~\ref{Q_6C_32} and \ref{Q_6 C_32 Decomp}, and the cycles  $C^1$, $C^2$, $B^1$, $B^2$, $Y^1$, and $Y^2$ are written explicitly at the end of the proof. The three partition sets in the decomposition are $\{C^1, C^2\}$, $\{B^1, Y^1\}$, and $\{B^2, Y^2\}$,

The graph $Q_6$ is represented in Figure~\ref{Q_6C_32} as the Cartesian product  $Q_4 \sq Q_2$. The cycles $C^1$ (dashed red) and $C^2$ (dotted blue) are shown in Figure~\ref{Q_6C_32}. Their vertex sets partition  $V(Q_6)$, so one partition set of our decomposition is $\{C^1, C^2\}$.

All edges not in $C^1$ or $C^2$ (gray, brown, and yellow edges, as well as edges not pictured) are represented in one of the isomorphic subgraphs $G_B$ and $G_Y$ of $Q_6$ in Figure~\ref{Q_6 C_32 Decomp} (top) and Figure~\ref{Q_6 C_32 Decomp} (middle) exactly once. Thus it remains to find a Hamiltonian decomposition of each of $G_B$ and $G_Y$, which is done in Figure~\ref{Q_6 C_32 Decomp} (bottom).  Let $B^1$ and $B^2$ be the cycles in the decomposition of $G_B$, and let $Y^1$ and $Y^2$ be the cycles in the decomposition of $G_Y$.
 
Since $V(G_B)$ and $V(G_Y)$ partition $V(Q_6)$, $\{B^1, Y^1\}$ and $\{B^2, Y^2\}$ are the remaining two partition sets.
 
The decomposition is summarized below.

The cycle $C^1$ is the  dashed red cycle in Figure~\ref{Q_6C_32}:
(011000, 011010, 011110, 001110, 001010, 101010, 101110, 111110, 111010, 111011, 111111, 101111, 101011, 001011, 001111, 011111, 011011, 011001, 011101, 001101, 001001, 101001, 101101, 111101, 111001, 111000, 111100, 101100, 101000, 001000, 001100, 011100, 011000)

The cycle $C^2$ is the dotted blue cycle in Figure~\ref{Q_6C_32}: (010100, 010110, 010010, 000010, 000110, 100110, 100010, 110010, 110110, 110111, 110011, 100011, 100111, 000111, 000011, 010011, 010111, 010101, 010001, 000001, 000101, 100101, 100001, 110001, 110101, 110100, 110000, 100000, 100100, 000100, 000000, 010000, 010100)

The cycles $B^1$ and $B^2$ are the two cycles in the Hamiltonian decomposition of the graph $G_B$ in Figure~\ref{Q_6 C_32 Decomp} (top).

$B^1$: (010110, 010111, 011111, 011101, 011100, 111100, 111101, 111111, 110111, 110101, 010101, 000101, 000111, 001111, 001101, 001100, 101100, 101101, 101111, 100111, 100101, 100100, 110100, 010100, 000100, 000110, 001110, 101110, 100110, 110110, 111110, 011110, 010110)

$B^2$: (010110, 110110, 110100, 111100, 111110, 111111, 011111, 011110, 011100, 010100, 010101, 011101, 111101, 110101, 100101, 101101, 001101, 000101, 000100, 001100, 001110, 001111, 101111, 101110, 101100, 100100, 100110, 100111, 110111, 010111, 000111, 000110, 010110)

The cycles $Y^1$ and $Y^2$ are the two cycles in the Hamiltonian decomposition of the graph $G_Y$ in Figure~\ref{Q_6 C_32 Decomp} (middle).

$Y^1$: (011010, 011011, 010011, 010001, 010000, 110000, 110001, 110011, 111011, 111001, 011001, 001001, 001011, 000011, 000001, 000000, 100000, 100001, 100011, 101011, 101001, 101000, 111000, 011000, 001000, 001010, 000010, 100010, 101010, 111010, 110010, 010010, 011010)

$Y^2$: (011010, 111010, 111000, 110000, 110010, 110011, 010011, 010010, 010000, 011000, 011001, 010001, 110001, 111001, 101001, 100001, 000001, 001001, 001000, 000000, 000010, 000011, 100011, 100010, 100000, 101000, 101010, 101011, 111011, 011011, 001011, 001010, 011010)
\end{proof}

Theorems~\ref{torus}, \ref{AS}, and \ref{ham} about Hamiltonian decompositions of Cartesian product graphs are due to Kotzig~\cite{kotzig} and Aubert and Schneider~\cite{aubert}, respectively.

\begin{theorem}[Kotzig, 1973]\label{torus}
The Cartesian product of two cycles has a Hamiltonian decomposition.
\end{theorem}

\begin{theorem}[Aubert and Schneider, 1982]\label{AS}
If $G$ is a 4-regular graph with a Hamiltonian decomposition and $C$ is a cycle, then $G \sq C$ has a Hamiltonian decomposition.
\end{theorem}

The fact that if $n$ is even, $Q_n$ has a Hamiltonian decomposition is a corollary of Theorem~\ref{AS} (see \cite{alspach}).  Note that a Hamiltonian decomposition is partitionable, with each partition set containing one cycle.  So we restate this result as follows.

\begin{theorem}\label{ham}
For all even $n$, $Q_n$ has a partitionable decomposition into Hamiltonian cycles.
\end{theorem}

We now state two propositions on decompositions of Cartesian products.  The first proposition is Proposition 8 in \cite{AOT21}.

\begin{proposition}\label{AOT}
Let $G$ and $H$ be graphs, with decompositions $G =G_1 \sqcup \dots \sqcup G_m$ and $H = H_1 \sqcup \dots \sqcup H_m$ where for $1 \le i \le m$, $V(G_i) = V(G)$ and $V(H_i) = V(H)$.  Then
\[G \sq H  = \bigsqcup_{i=1}^m( G_i \sq H_i).\]
\end{proposition}

\begin{proposition}\label{decomp}
Let $G$ be a graph with a partitionable decomposition $G =G_1 \sqcup \dots \sqcup G_m$ with one partition set $\F = \{G_1, \ldots, G_m\}$, and let $H$ be any graph.
Then 
\[G \sq H = \bigsqcup_{i=1}^m (G_i \sq H)\]
is a partitionable decomposition of $G \sq H$ with one partition set.

\end{proposition}

\begin{proof}
We first show that $E(G \sq H) = \bigsqcup_{i=1}^m E(G_i \sq H)$. Fix an edge $e \in E(G \sq H)$ and consider two cases, depending whether $e$ is a horizontal or vertical edge.

If $e$ is horizontal, $e = (u, v)(u', v)$, where $uu' \in E(G), v \in V(H)$. Then because $G =G_1 \sqcup \dots \sqcup G_m$, we know $uu' \in E(G_j)$ for exactly one $j \in \{1, \ldots, m\}$, and therefore $e \in E(G_i \sq H)$ if and only if $i=j$.

Otherwise $e$ is vertical, so $e = (u, v)(u, v')$, where $vv' \in E(H), u \in V(G)$. Then since $\F = \{G_1, \ldots, G_m\}$ is a partition set, we know $u \in V(G_j)$ for exactly one $j \in \{1, \ldots, m\}$ and therefore $e \in E(G_i \sq H)$ if and only if $i=j$.

It remains to show that $\{G_1 \sq H, \ldots, G_m \sq H\}$ is a partition set, i.e. their vertex sets partition  $V(G \sq H)$. 
Fix a vertex $(u, v) \in V(G \sq H)$. Because $\F = \{G_1, \ldots , G_m\}$ is a partition set for $G$, we know $u \in V(G_j)$ for exactly one $j \in \{1, \ldots, m\}$. Therefore $(u, v) \in V(G_i \sq H)$ if and only if $i=j$.
\end{proof}

In Lemma~\ref{hhlemma} and Proposition~\ref{cart} we will use the following observation, stating roughly that a partitionable decomposition of a partitionable decomposition into spanning subgraphs gives a partitionable decomposition.

\begin{proposition}\label{ddd}
Suppose $G$ has a partitionable decomposition $G = G_1 \sqcup \dots \sqcup G_m$, and for $1\le i \le m$, $V(G_i) = V(G)$. Suppose for $1\le i \le m$, $G_i$ has a partitionable decomposition into copies of a graph $H$, with partition sets $\F_{i,1}, \ldots, \F_{i,k_i}$. Then $G$ has a partitionable decomposition into copies of $H$, with partition sets $\F_{i,j}$ where $1 \le i \le m$ and $1 \le j \le k_i$.
\end{proposition}

\begin{lemma}\label{hhlemma}
Suppose the graph $G$ has a partitionable decomposition into cycles of length $\ell$, with $m$ partition sets, where $m \ge 1$. Suppose the graph $H$ has a Hamiltonian decomposition with between $m$ and $2m$ cycles.  Then $G \sq H$ has a partitionable decomposition into cycles of length $\ell |V(H)|$.
\end{lemma}

\begin{proof}
Suppose $G$ has a partitionable decomposition into cycles of length $\ell$ with $m$ partition sets $\F_1, \dots, \F_m$.  For $1 \le i \le m$, let $G_i = \cup_{C \in \F_i}C$. Since $\F_i$ is a partition set, for $1 \le i\le m$, $V(G_i) = V(G)$.
Suppose $H$ has a  Hamiltonian decomposition with $m+n$ cycles, where $0 \le n \le m$. Denote these Hamiltonian cycles by $H^1, \dots, H^{m+n}$. Then $H$ has a decomposition $H = H_1 \sqcup \dots \sqcup H_m$, where for $1\le i \le n$, $H_i=H^i \cup H^{m+i}$, and for $n+1 \le i \le m$, $H_i = H^{i}$. 

Since $H_i$ is either a Hamiltonian cycle or a union of two Hamiltonian cycles, for $1\le i \le m$, $V(H_i) = V(H)$. Thus by Proposition~\ref{AOT}, $G \sq H = \sqcup_{i=1}^m G_i \sq H_i$.

For $1 \le i \le m$, since $V(G_i) = V(G)$ and $V(H_i) = V(H)$, $V(G_i \sq H_i) = V(G \sq H)$.  Thus by Proposition~\ref{ddd}, it remains to show that for $1 \le i \le m$, $G_i \sq H_i$ has a partitionable decomposition into cycles of length $\ell|V(H)|$.   
By Proposition~\ref{decomp}, letting $G = G_i$, $\F = \F_i$, and $H = H_i$,
\[G_i \sq H_i = \bigsqcup_{C \in \F_i} (C \sq H_i),\]
and this decomposition is partitionable with one partition set.

There are now two cases for the decomposition of each Cartesian product $C \sq H_i$, depending if $H_i$ is the union of one or two cycles.

Case 1:  If $n+ 1 \le i \le m$, then $H_i$ is a Hamiltonian cycle in $H$.  Each graph in the union $\bigsqcup_{C \in \F_i} (C \sq H_i)$ is a Cartesian product of two cycles of lengths $\ell$ and $|V(H)|$, which by Theorem~\ref{torus} has a decomposition into two Hamiltonian cycles, each of length $\ell|V(H)|$.  Since the decomposition of $G_i \sq H_i$ has one partition set, if we color the Hamiltonian cycles from each torus red and blue, the set of red cycles forms one partition set, while the set of blue cycles forms another.

Case 2: If $1 \le i \le n$, $H_i$ is a union of two Hamiltonian cycles.  Theorem~\ref{AS} implies each graph in the union $\bigsqcup_{C \in \F_i} (C \sq H_i)$ can be decomposed into three Hamiltonian cycles, each of length $\ell|V(H)|$.  Since the decomposition of $G_i \sq H_i$ has one partition set, if we color the Hamiltonian cycles from each $C \sq H_i$ red, blue, and green the set of red cycles forms one partition set, the set of blue cycles forms another, and the set of green cycles forms a third.
\end{proof}

\begin{proposition}\label{halfham}
If $n$ is even and $n \ge 4$, then $Q_n$ has a partitionable decomposition into cycles of length $2^{n-1}$.
\end{proposition}

\begin{proof}
The proof is by induction on even $n$. The two base cases $n = 4$ and $n = 6$ are given by Propositions~\ref{q4} and \ref{q6}.

For the inductive step, fix even $n \geq 6$ and assume $Q_k$ has a partitionable decomposition into cycles of length $2^{k-1}$ for all even $k$ where $4 \le k \leq n$. We show that $Q_{n+2}$ has a partitionable decomposition into cycles of length $2^{(n+2)-1}$ by considering two cases depending on whether $n + 2$ is a multiple of 4.   Let $Q_{n+2} = Q_{4m + 2i} = Q_{2m} \sq Q_{2m + 2i}$, for some natural number $m$, where $i \in \{0, 1\}$. Since $n+2 \ge 8$, $2m$ must be at least 4, so by the induction hypothesis $Q_{2m}$ has a partitionable decomposition into cycles of length $2^{2m-1}$, with $m$ partition sets, each containing two cycles. By Theorem~\ref{ham} we know that $Q_{2m + 2i}$ has a partionable (Hamiltonian) decomposition into cycles of length $2^{2m + 2i}$, with $m + i$ partition sets (each containing one cycle). Since $m \leq m + i \leq 2m$, letting $G = Q_{2m}$, $\ell =  2^{2m-1}$, and $H = Q_{2m+2i}$, Lemma~\ref{hhlemma} implies $Q_{n+2} = Q_{2m} \sq Q_{2m + 2i}$ has a partitionable decomposition into cycles of length $2^{2m-1} |V(Q_{2m+2i})| = 2^{2m-1} \cdot 2^{2m + 2i} = 2^{4m +2i -  1} = 2^{(n+2)-1}$.
\end{proof}

We now prove Lemma~\ref{8ell} and Proposition~\ref{cart}, which are the results needed for the inductive argument in the proof of Theorem~\ref{partitionable}.

\begin{lemma}\label{8ell}
If $\ell \ge 1$, $n \ge 1$, and $n$ divides $4\ell$, then  $C_{4\ell} \sq C_4$ has a partitionable decomposition into cycles of length $4n$ with two partition sets.
\end{lemma}

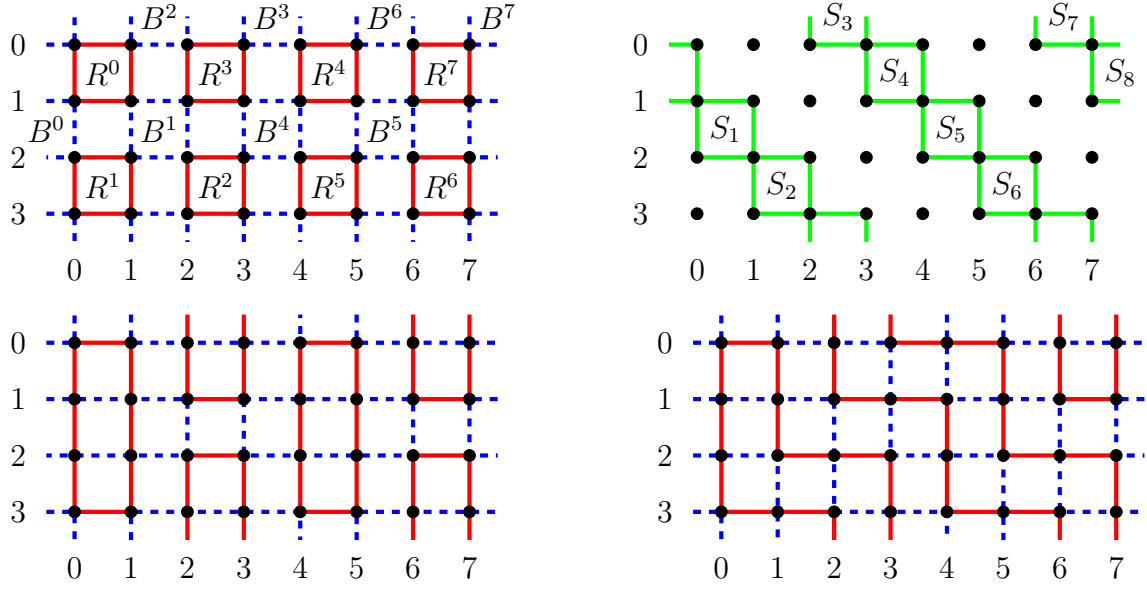
\begin{figure}
\begin{tikzpicture}[scale = .75]

\node at (0,-1){0};
\node at (1,-1){1};
\node at (2,-1){2};
\node at (3,-1){3};
\node at (4,-1){4};
\node at (5,-1){5};
\node at (6,-1){6};
\node at (7,-1){7};

\node at (-1,3){0};
\node at (-1,2){1};
\node at (-1,1){2};
\node at (-1,0){3};

\draw[ultra thick, red]
(0, 0) -- (0, 1) -- (1,1) -- (1, 0) -- (0,0);
\draw[ultra thick, red]
(0, 2) -- (0, 3) -- (1,3) -- (1, 2) -- (0,2);
\draw[ultra thick, red]
(2, 0) -- (2, 1) -- (3,1) -- (3,0) -- (2,0);
\draw[ultra thick, red]
(2, 2) -- (2, 3) -- (3,3) -- (3,2) -- (2,2);
\draw[ultra thick, red]
(4, 0) -- (4, 1) -- (5,1) -- (5, 0) -- (4,0);
\draw[ultra thick, red]
(4, 2) -- (4, 3) -- (5,3) -- (5, 2) -- (4,2);
\draw[ultra thick, red]
(6, 0) -- (6, 1) -- (7,1) -- (7,0) -- (6,0);
\draw[ultra thick, red]
(6, 2) -- (6, 3) -- (7,3) -- (7,2) -- (6,2);

\draw[ultra thick, blue, dashed]
(-.5, 2) -- (0, 2) -- (0, 1) -- (-.5, 1);
\draw[ultra thick, blue, dashed]
(7.5, 2) -- (7, 2) -- (7,1) -- (7.5,1);
\draw[ultra thick, blue, dashed]
(1, -.5) -- (1, 0) -- (2, 0) -- (2, -.5);
\draw[ultra thick, blue, dashed]
(1, 3.5) -- (1, 3) -- (2, 3) -- (2, 3.5);
\draw[ultra thick, blue, dashed]
(3, -.5) -- (3, 0) -- (4, 0) -- (4, -.5);
\draw[ultra thick, blue, dashed]
(3, 3.5) -- (3, 3) -- (4, 3) -- (4, 3.5);
\draw[ultra thick, blue, dashed]
(1, 2) -- (2,2) -- (2,1) -- (1,1) -- (1,2);
\draw[ultra thick, blue, dashed]
(3, 2) -- (4,2) -- (4,1) -- (3,1) -- (3,2);
\draw[ultra thick, blue, dashed]
(5, 2) -- (6,2) -- (6,1) -- (5,1) -- (5,2);
\draw[ultra thick, blue, dashed]
(5, 3.5) -- (5, 3) -- (6, 3) -- (6, 3.5);
\draw[ultra thick, blue, dashed]
(5, -.5) -- (5, 0) -- (6, 0) -- (6,-.5);
\draw[ultra thick, blue, dashed]
(-.5, 0) -- (0, 0) -- (0, -.5);
\draw[ultra thick, blue, dashed]
(-.5, 3) -- (0, 3) -- (0, 3.5);
\draw[ultra thick, blue, dashed]
(7.5, 0) -- (7, 0) -- (7, -.5);
\draw[ultra thick, blue, dashed]
(7, 3.5) -- (7, 3) -- (7.5, 3);

\node at (.5, 2.5){$R^0$};
\node at (.5, .5){$R^1$};
\node at (-.5, 1.5){$B^0$};
\node at (1.5, 1.5){$B^1$};
\node at (2.5, 2.5){$R^3$};
\node at (2.5, .5){$R^2$};
\node at (1.5, 3.5){$B^2$};
\node at (3.5, 3.5){$B^3$};
\node at (4.5, 2.5){$R^4$};
\node at (4.5, .5){$R^5$};
\node at (3.5, 1.5){$B^4$};
\node at (5.5, 1.5){$B^5$};
\node at (6.5, 2.5){$R^7$};
\node at (6.5, .5){$R^6$};
\node at (5.5, 3.5){$B^6$};
\node at (7.5, 3.5){$B^7$};

\foreach \x in {0, ..., 7}{
    \foreach \y in {0,...,3}{
    \node[style={circle, draw, fill=black},scale=.4] at (\x,\y){};}};

\end{tikzpicture}
\hfill
\begin{tikzpicture}[scale = .75]

\node at (0,-1){0};
\node at (1,-1){1};
\node at (2,-1){2};
\node at (3,-1){3};
\node at (4,-1){4};
\node at (5,-1){5};
\node at (6,-1){6};
\node at (7,-1){7};

\node at (-1,3){0};
\node at (-1,2){1};
\node at (-1,1){2};
\node at (-1,0){3};

\draw[ultra thick, green] (0,2) -- (1,2) -- (1,1) -- (0,1) -- (0,2);

\draw[ultra thick, green] (1,0) -- (1,1) -- (2,1) -- (2,0) -- (1,0);

\draw[ultra thick, green] (2,-.5) -- (2,0) -- (3,0) -- (3,-.5);

\draw[ultra thick, green] (2,3.5) -- (2,3) -- (3,3) -- (3,3.5);

\draw[ultra thick, green] (3,2) -- (3,3) -- (4,3) -- (4,2) -- (3,2);

\draw[ultra thick, green] (4,2) -- (5,2) -- (5,1) -- (4,1) -- (4,2);

\draw[ultra thick, green] (5,0) -- (5,1) -- (6,1) -- (6,0) -- (5,0);

\draw[ultra thick, green] (6,-.5) -- (6,0) -- (7,0) -- (7,-.5);

\draw[ultra thick, green] (6,3.5) -- (6,3) -- (7,3) -- (7,3.5);

\draw[ultra thick, green] (7.5, 3) -- (7,3) -- (7,2) -- (7.5,2);

\draw[ultra thick, green] (-.5, 3) -- (0,3) -- (0,2) -- (-.5,2);

\foreach \x in {0, ..., 7}{
    \foreach \y in {0,...,3}{
    \node[style={circle, draw, fill=black},scale=.4] at (\x,\y){};}};
    
\node at (.5, 1.5){$S_1$};
\node at (1.5, .5){$S_2$};
\node at (2.5, 3.5){$S_3$};
\node at (3.5, 2.5){$S_4$};
\node at (4.5, 1.5){$S_5$};
\node at (5.5, .5){$S_6$};
\node at (6.5, 3.5){$S_7$};
\node at (7.5, 2.5){$S_8$};

\end{tikzpicture}

\vspace{.1in}

\begin{tikzpicture}[scale = .75]

\node at (0,-1){0};
\node at (1,-1){1};
\node at (2,-1){2};
\node at (3,-1){3};
\node at (4,-1){4};
\node at (5,-1){5};
\node at (6,-1){6};
\node at (7,-1){7};

\node at (-1,3){0};
\node at (-1,2){1};
\node at (-1,1){2};
\node at (-1,0){3};

\draw[ultra thick, red]
(0, 3) -- (1,3) -- (1,0) -- (0,0) -- (0,3);

\draw[ultra thick, red]
(4, 3) -- (5,3) -- (5,0) -- (4,0) -- (4,3);

\draw[ultra thick, red]
(2, -.5) -- (2, 1) -- (3,1) -- (3, -.5);

\draw[ultra thick, red]
(2, 3.5) -- (2, 2) -- (3,2) -- (3, 3.5);

\draw[ultra thick, red]
(6, -.5) -- (6, 1) -- (7,1) -- (7, -.5);

\draw[ultra thick, red]
(6, 3.5) -- (6, 2) -- (7,2) -- (7, 3.5);

\draw[ultra thick, blue, dashed]
(-.5, 2) -- (2, 2) -- (2, 1) -- (-.5, 1);

\draw[ultra thick, blue, dashed] (7.5, 2) -- (7,2) -- (7, 1) -- (7.5, 1);

\draw[ultra thick, blue, dashed] (3, 2) -- (3,1) -- (6, 1) -- (6,2) -- (3,2);

\draw[ultra thick, blue, dashed] (1, 3.5) -- (1,3) -- (4, 3) -- (4,3.5);

\draw[ultra thick, blue, dashed] (1, -.5) -- (1,0) -- (4, 0) -- (4,-.5);

\draw[ultra thick, blue, dashed] (5, 3.5) -- (5,3) -- (7.5, 3);

\draw[ultra thick, blue, dashed] (5, -.5) -- (5, 0) -- (7.5, 0);

\draw[ultra thick, blue, dashed] (-.5, 3) -- (0,3) -- (0, 3.5);

\draw[ultra thick, blue, dashed] (-.5, 0) -- (0,0) -- (0, -.5);

\foreach \x in {0, ..., 7}{
    \foreach \y in {0,...,3}{
    \node[style={circle, draw, fill=black},scale=.4] at (\x,\y){};}};
\end{tikzpicture}
\hfill
\begin{tikzpicture}[scale = .75]

\node at (0,-1){0};
\node at (1,-1){1};
\node at (2,-1){2};
\node at (3,-1){3};
\node at (4,-1){4};
\node at (5,-1){5};
\node at (6,-1){6};
\node at (7,-1){7};

\node at (-1,3){0};
\node at (-1,2){1};
\node at (-1,1){2};
\node at (-1,0){3};

\draw[ultra thick, red]
(2, -.5) -- (2,0) -- (0,0) -- (0,3) -- (1,3) -- (1,1) -- (3,1) -- (3,-.5);

\draw[ultra thick, red]
(2, 3.5) -- (2,2) -- (4, 2) -- (4,0) -- (6,0) -- (6,-.5);

\draw[ultra thick, red]
(3,3.5) -- (3,3) -- (5,3) -- (5,1) -- (7,1) -- (7,-.5);

\draw[ultra thick, red]
(6,3.5) -- (6,2) -- (7,2) -- (7,3.5);

\draw[ultra thick, blue, dashed]
(-.5, 2) -- (2,2) -- (2,0) -- (4,0) -- (4,-.5);

\draw[ultra thick, blue, dashed]
(-.5, 1) -- (1,1) -- (1, -.5);

\draw[ultra thick, blue, dashed]
(-.5, 0) -- (0,0) -- (0,-.5);

\draw[ultra thick, blue, dashed]
(-.5, 3) -- (0,3) -- (0,3.5);

\draw[ultra thick, blue, dashed]
(5, 3.5) -- (5,3) -- (7.5, 3);

\draw[ultra thick, blue, dashed]
(4, 3.5) -- (4, 2) -- (6,2) -- (6,0) -- (7.5, 0);

\draw[ultra thick, blue, dashed]
(1, 3.5) -- (1,3) -- (3,3) -- (3,1) -- (5,1) -- (5,-.5);

\draw[ultra thick, blue, dashed]
(7.5, 2) -- (7,2) -- (7,1) -- (7.5,1);

\foreach \x in {0, ..., 7}{
    \foreach \y in {0,...,3}{
    \node[style={circle, draw, fill=black},scale=.4] at (\x,\y){};}};
\end{tikzpicture}

\captionof{figure}{Top left: The partitionable decomposition of $C_{8} \sq C_4$ into cycles of length four from the proof of Lemma~\ref{8ell} ($\ell=2$). Top right: Potential recoloring locations $S_1, \ldots, S_8$. Bottom left: A partitionable decomposition of $C_{8} \sq C_4$ into cycles of length 8. Since $n = 2$ the recolored cycles are $L = \{S_1, S_3, S_5, S_7\}$. Bottom right: A partitionable decomposition of $C_{8} \sq C_4$ into cycles of length 32. Since $n = 8$, the recolored cycles are $L = \{S_1, S_2, S_3, S_4, S_5, S_6, S_7\}$.}
\label{C8C4n1}
\end{figure}

Before proving the lemma, note that Figures~\ref{C8C4n1} (bottom left and bottom right) and \ref{C_24 sq C_4, recoloring & initial} (bottom middle and bottom) contain examples of decompositions of $C_{4\ell} \sq C_4$ into cycles of length $4n$ where $(n, \ell)=(2,2)$, $(n, \ell)=(8,2)$, $(n, \ell)=(6,6)$, and $(n, \ell)=(24,6)$, respectively.

\begin{proof}
Since $C_{4\ell} \sq C_4$ is a 4-regular graph with $16\ell$ vertices, a partitionable decomposition of $C_{4\ell} \sq C_4$ into cycles of length $4n$ will have 2 partition sets, each of which contain $4\ell/n$ cycles.
We start with a partitionable decomposition of $C_{4\ell} \sq C_4$ into cycles of length four, then show how to modify this decomposition to obtain a partitionable decomposition of $C_{4\ell} \sq C_4$ into cycles of length $4n$.

Let $V(C_{4\ell}) = \{0,1, \ldots, 4\ell - 1\}$ and $V(C_4) = \{0,1,2,3\}$, so $V(C_{4\ell} \sq C_4) = \{0,1,\ldots, 4\ell - 1\} \times \{0,1,2,3\}$. For $0 \le k \le 4\ell-1$, we describe ``red'' and ``blue'' cycles $R^k$ and $B^k$ of length four. Let $R^0 = ((0,0), (0, 1), (1,1), (1,0), (0,0))$, $R^1 = ((0,2), (0, 3), (1,3), (1,2), (0,2))$, $B^0 = ((4\ell - 1, 1), (4\ell - 1, 2), (0,2),(0, 1) , (4\ell - 1, 1))$, and $B^1 = ((1, 1), (1,2), (2,2), (2, 1), (1, 1))$. For $2 \leq k \leq 4 \ell - 1$, given a cycle $C^{k-2}$ (where $C$ could be $R$ or $B$), let $C^{k}$ be the cycle such that $V(C^{k}) = \{(i+2 \text{ (mod $4\ell$)}, j+2 \text{ (mod 4))} : (i,j) \in V(C^{k-2})\}$. 
Figure~\ref{C8C4n1} (top left) shows this decomposition on $C_8 \sq C_4$ (where $\ell = 2$), while Figure ~\ref{C_24 sq C_4, recoloring & initial} (top) shows this decomposition on $C_{24} \sq C_4$ (where $\ell = 6$).

We verify that $\{R^0, \ldots, R^{4\ell-1}, B^0, \ldots, B^{4\ell-1}\}$ forms a decomposition of $C_{4\ell} \sq C_4$: By inspection, we can verify that for $0 \le i, j \le 3$, the vertex $(i,j)$ is incident with two red edges and two blue edges (see for example Figure~\ref{C8C4n1}).  For  $0 \le i \le 4\ell-1$ and $0 \le j \le 3$, the vertex $(i,j)$ is incident with edges of the same color as $(i',j)$, where $i \equiv i' \mod{4}$.  Thus every edge is in exactly one cycle (i.e. is colored either red or blue), and we have a decomposition. The same observation shows that each vertex is in exactly one red cycle and one blue cycle, so $\{R^0, \ldots, R^{4\ell-1}\}$ and $\{B^0, \ldots, B^{4\ell-1}\}$ are the partition sets of this partitionable decomposition.

Now we show how to transform our partitionable decomposition of cycles of length four into a partitionable decomposition of cycles of length $4n$. This can be accomplished by recoloring some of the edges.

We call the following recoloring a cycle combination operation: Suppose we are given vertex-disjoint red cycles $R$ and $R'$ where $v$ and $w$  are consecutive vertices in $R$ and $v'$ and $w'$ are consecutive vertices in $R'$.  We also have  vertex-disjoint blue cycles $B$ and $B'$ where $v$ and $v'$  are consecutive vertices in $B$ and $w$ and $w'$ are consecutive vertices in $B'$. Recoloring $vw$ and $v'w'$ to be blue and $vv'$ and $ww'$ to be red yields a single red cycle whose length is the sum of the lengths of $R$ and $R'$ and contains each vertex from $R$ and $R'$ exactly once and a single blue cycle  whose length is the sum of the lengths of $B$ and $B'$ and contains each vertex from $B$ and $B'$ exactly once.  I.e. recoloring the 4-cycle $(v, w, w', v', v)$ combines the red cycles and combines the blue cycles.

We now identify specific 4-cycles $S_1, \ldots, S_{4\ell}$ in $C_{4\ell} \sq C_4$ where the cycle combination operation can be applied. Let $S_1$ = $((0,1), (0,2), (1,2), (1,1), (0,1))$ and for $2 \leq k \leq 4 \ell$, let $S_{k}$ be such that $V(S_{k}) = \{(i+1 \text{ (mod $4\ell$)}, j+1 \text{ (mod 4)}) : (i,j) \in V(S_{k-1})\}$. Figure~\ref{C8C4n1} (top right) illustrates the eight cycles $S_1, \ldots, S_8$ in $C_8 \sq C_4$, while Figure~\ref{C_24 sq C_4, recoloring & initial} (top middle) illustrates the 24 cycles $S_1, \ldots,  S_{24}$ in $C_{24} \sq C_4$.  Note that recoloring $S_i$ applies the cycle combination operation to the red cycles $R^{i-1}$ and $R^i$ and the blue cycles $B^{i-1}$ and $B^{i}$, where the superscripts are interpreted mod $4\ell$.

Thus if $n$ divides $4\ell$ we can create a partitionable decomposition of $C_{4\ell} \sq C_4$ into cycles of length $4n$ by recoloring sets of $n-1$ consecutive recoloring locations. Let $L = \{S_i: n \text{ does not divide } i\}$. 
The set $L$ contains disjoint sets of $n-1$ consecutive recoloring locations, so recoloring all of the cycles $S_i$ where $S_i \in L$ gives red and blue cycles of length $4n$. Since applying the cycle combination operation leaves vertices in the recolored red cycles the same as in the original red cycles, and the same is true for the blue cycles, the red cycles of length $4n$ form a partition set, as do the blue cycles of length $4n$.  
\end{proof}

\begin{figure}

\begin{tikzpicture}[scale = 1.2]
\node at (0,-.5){0};
\node at (1,-.5){2};
\node at (2,-.5){4};
\node at (3,-.5){6};
\node at (4,-.5){8};
\node at (5,-.5){10};
\node at (6,-.5){12};
\node at (7,-.5){14};
\node at (8,-.5){16};
\node at (9,-.5){18};
\node at (10,-.5){20};
\node at (11,-.5){22};

\node at (-.5,1.5){0};
\node at (-.5,1){1};
\node at (-.5,.5){2};
\node at (-.5,0){3};

\foreach \x in {0, ..., 11}{
\draw[ultra thick, red]
(\x,1.5) -- (\x +.5,1.5) -- (\x + .5,1) -- (\x,1) -- (\x,1.5);
\draw[ultra thick, red]
(\x,.5) -- (\x +.5,.5) -- (\x + .5,0) -- (\x,0) -- (\x,.5);}

\draw[ultra thick, blue, dashed]
(-.25, 1) -- (0,1) -- (0,.5) -- (-.25, .5);

\draw[ultra thick, blue, dashed]
(11.75, 1) -- (11.5,1) -- (11.5,.5) -- (11.75, .5);

\foreach \x in {0, ..., 10}{
\draw[ultra thick, blue, dashed]
(\x + .5,1) -- (\x +1,1) -- (\x + 1,.5) -- (\x + .5,.5) -- (\x + .5 ,1);

\draw[ultra thick, blue, dashed](\x + .5, 1.75) -- (\x + .5, 1.5) -- (\x + 1, 1.5) -- (\x + 1, 1.75);

\draw[ultra thick, blue, dashed](\x + .5, -.25) -- (\x + .5, 0) -- (\x + 1, 0) -- (\x + 1, -.25);
}

\node at (.25,1.25){$R^{0}$};
\node at (1.25,1.25){$R^{3}$};
\node at (2.25,1.25){$R^{4}$};
\node at (3.25,1.25){$R^{7}$};
\node at (4.25,1.25){$R^{8}$};
\node at (5.25,1.25){$R^{11}$};
\node at (6.25,1.25){$R^{12}$};
\node at (7.25,1.25){$R^{15}$};
\node at (8.25,1.25){$R^{16}$};
\node at (9.25,1.25){$R^{19}$};
\node at (10.25,1.25){$R^{20}$};
\node at (11.25,1.25){$R^{23}$};

\node at (.25,.25){$R^{1}$};
\node at (1.25,.25){$R^{2}$};
\node at (2.25,.25){$R^{5}$};
\node at (3.25,.25){$R^{6}$};
\node at (4.25,.25){$R^{9}$};
\node at (5.25,.25){$R^{10}$};
\node at (6.25,.25){$R^{13}$};
\node at (7.25,.25){$R^{14}$};
\node at (8.25,.25){$R^{17}$};
\node at (9.25,.25){$R^{18}$};
\node at (10.25,.25){$R^{21}$};
\node at (11.25,.25){$R^{22}$};

\node at 
(-.25, .75){$B^{0}$};
\node at 
(.75, .75){$B^{1}$};
\node at 
(1.75, .75){$B^{4}$};
\node at 
(2.75, .75){$B^{5}$};
\node at 
(3.75, .75){$B^{8}$};
\node at 
(4.75, .75){$B^{9}$};
\node at 
(5.75, .75){$B^{12}$};
\node at 
(6.75, .75){$B^{13}$};
\node at 
(7.75, .75){$B^{16}$};
\node at 
(8.75, .75){$B^{17}$};
\node at 
(9.75, .75){$B^{20}$};
\node at 
(10.75, .75){$B^{21}$};

\node at 
(.75, 1.75){$B^{2}$};
\node at 
(1.75, 1.75){$B^{3}$};
\node at 
(2.75, 1.75){$B^{6}$};
\node at 
(3.75, 1.75){$B^{7}$};
\node at 
(4.75, 1.75){$B^{10}$};
\node at 
(5.75, 1.75){$B^{11}$};
\node at 
(6.75, 1.75){$B^{14}$};
\node at 
(7.75, 1.75){$B^{15}$};
\node at 
(8.75, 1.75){$B^{18}$};
\node at 
(9.75, 1.75){$B^{19}$};
\node at 
(10.75, 1.75){$B^{22}$};
\node at 
(11.75, 1.75){$B^{23}$};

\draw[ultra thick, blue, dashed]
(-.25, 1.5) -- (0,1.5) -- (0,1.75);

\draw[ultra thick, blue, dashed]
(-.25, 0) -- (0,0) -- (0,-.25);

\draw[ultra thick, blue, dashed]
(11.75, 1.5) -- (11.5,1.5) -- (11.5,1.75);

\draw[ultra thick, blue, dashed]
(11.75, 0) -- (11.5,0) -- (11.5,-.25);

;

\foreach \x in {0, .5, 1, 1.5, 2, 2.5, 3, 3.5}{
    \foreach \y in {0,.5, 1,1.5}{
    \node[style={circle, draw, fill=black},scale=.4] at (\x,\y){};}};

\foreach \x in {4, 4.5, 5, 5.5, 6, 6.5, 7, 7.5}{
    \foreach \y in {0,.5,1,1.5}{
    \node[style={circle, draw, fill=black},scale=.4] at (\x,\y){};}};

\foreach \x in {8, 8.5, 9, 9.5, 10, 10.5, 11, 11.5}{
    \foreach \y in {0,.5,1,1.5}{
    \node[style={circle, draw, fill=black},scale=.4] at (\x,\y){};}};
\end{tikzpicture}

\vspace{.1in}

\begin{tikzpicture}[scale = 1.2]

\node at (0,-.5){0};
\node at (1,-.5){2};
\node at (2,-.5){4};
\node at (3,-.5){6};
\node at (4,-.5){8};
\node at (5,-.5){10};
\node at (6,-.5){12};
\node at (7,-.5){14};
\node at (8,-.5){16};
\node at (9,-.5){18};
\node at (10,-.5){20};
\node at (11,-.5){22};

\node at (-.5,1.5){0};
\node at (-.5,1){1};
\node at (-.5,.5){2};
\node at (-.5,0){3};

\foreach \x in {0, 2, 4, 6, 8}{
\draw[ultra thick, green]
(\x, .5) -- (\x, 1) -- (\x + .5, 1) -- (\x + .5, .5) -- (\x, .5);

\draw[ultra thick, green]
(\x + .5, .5) -- (\x + 1, .5) -- (\x + 1, 0) -- (\x + .5, 0) -- (\x + .5, .5);

\draw[ultra thick, green]
(\x + 1, -.25) -- (\x + 1, 0) -- (\x + 1.5, 0) -- (\x + 1.5, -.25);

\draw[ultra thick, green]
(\x + 1, 1.75) -- (\x + 1, 1.5) -- (\x + 1.5, 1.5) -- (\x + 1.5, 1.75);

\draw[ultra thick, green]
(\x + 1.5, 1.5) -- (\x + 2, 1.5) -- (\x + 2, 1) -- (\x + 1.5, 1) -- (\x + 1.5, 1.5);

}

\draw[ultra thick, green]
(10, .5) -- (10, 1) -- (10 + .5, 1) -- (10 + .5, .5) -- (10, .5);

\draw[ultra thick, green]
(10 + .5, .5) -- (10 + 1, .5) -- (10 + 1, 0) -- (10 + .5, 0) -- (10 + .5, .5);

\draw[ultra thick, green]
(10 + 1, -.25) -- (10 + 1, 0) -- (10 + 1.5, 0) -- (10 + 1.5, -.25);

\draw[ultra thick, green]
(10 + 1, 1.75) -- (10 + 1, 1.5) -- (10 + 1.5, 1.5) -- (10 + 1.5, 1.75);

\draw[ultra thick, green]
(-.25, 1.5) -- (0,1.5) -- (0,1) -- (-.25, 1);

\draw[ultra thick, green]
(11.75, 1.5) -- (11.5, 1.5) -- (11.5, 1) -- (11.75, 1);

\node at (.25, .75){$S_1$};
\node at (.75, .25){$S_2$};
\node at (1.25, 1.75){$S_3$};
\node at (1.75, 1.25){$S_4$};
\node at (2.25, .75){$S_5$};
\node at (2.75, .25){$S_6$};

\node at (3.25, 1.75){$S_7$};
\node at (3.75, 1.25){$S_8$};
\node at (4.25, .75){$S_9$};
\node at (4.75, .25){$S_{10}$};

\node at (5.25, 1.75){$S_{11}$};
\node at (5.75, 1.25){$S_{12}$};
\node at (6.25, .75){$S_{13}$};
\node at (6.75, .25){$S_{14}$};

\node at (7.25, 1.75){$S_{15}$};
\node at (7.75, 1.25){$S_{16}$};
\node at (8.25, .75){$S_{17}$};
\node at (8.75, .25){$S_{18}$};

\node at (9.25, 1.75){$S_{19}$};
\node at (9.75, 1.25){$S_{20}$};
\node at (10.25, .75){$S_{21}$};
\node at (10.75, .25){$S_{22}$};

\node at (11.25, 1.75){$S_{23}$};

\node at (11.75, 1.25){$S_{24}$};

\foreach \x in {0, .5, 1, 1.5, 2, 2.5, 3, 3.5}{
    \foreach \y in {0,.5, 1,1.5}{
    \node[style={circle, draw, fill=black},scale=.4] at (\x,\y){};}};

\foreach \x in {4, 4.5, 5, 5.5, 6, 6.5, 7, 7.5}{
    \foreach \y in {0,.5,1,1.5}{
    \node[style={circle, draw, fill=black},scale=.4] at (\x,\y){};}};

\foreach \x in {8, 8.5, 9, 9.5, 10, 10.5, 11, 11.5}{
    \foreach \y in {0,.5,1,1.5}{
    \node[style={circle, draw, fill=black},scale=.4] at (\x,\y){};}};

\end{tikzpicture}

\vspace{.1in}

\begin{tikzpicture}[scale = 1.2]

\node at (-.5,1.5){0};
\node at (-.5,1){1};
\node at (-.5,.5){2};
\node at (-.5,0){3};

\node at (0,-.5){0};
\node at (1,-.5){2};
\node at (2,-.5){4};
\node at (3,-.5){6};
\node at (4,-.5){8};
\node at (5,-.5){10};
\node at (6,-.5){12};
\node at (7,-.5){14};
\node at (8,-.5){16};
\node at (9,-.5){18};
\node at (10,-.5){20};
\node at (11,-.5){22};

\draw[ultra thick, red]
(1, -.25) -- (1, 0) -- (0,0) -- (0, 1.5) -- (.5, 1.5) -- (.5, .5) -- (1.5, .5) -- (1.5, -.25);

\draw[ultra thick, red]
(1, 1.75) -- (1,1) -- (2, 1) -- (2, 0) -- (2.5,0) -- (2.5, 1.5) -- (1.5,1.5) -- (1.5, 1.75);

\draw[ultra thick, red]
(3, -.25) -- (3, .5) -- (3.5,.5) -- (3.5, -.25);

\draw[ultra thick, red]
(3, 1.75) -- (3, 1) -- (4, 1) -- (4,0) -- (5,0) -- (5, -.25);

\draw[ultra thick, red]
(3.5, 1.75) -- (3.5, 1.5) -- (4.5, 1.5) -- (4.5, .5) -- (5.5, .5) -- (5.5, -.25);

\draw[ultra thick, red]
(5, 1.75) -- (5, 1) -- (5.5, 1) -- (5.5, 1.75);

\draw[ultra thick, red]
(7, -.25) -- (7, 0) -- (6,0) -- (6, 1.5) -- (6.5, 1.5) -- (6.5, .5) -- (7.5, .5) -- (7.5, -.25);

\draw[ultra thick, red]
(7, 1.75) -- (7,1) -- (8, 1) -- (8, 0) -- (8.5,0) -- (8.5, 1.5) -- (7.5,1.5) -- (7.5, 1.75);

\draw[ultra thick, red]
(9, -.25) -- (9, .5) -- (9.5,.5) -- (9.5, -.25);

\draw[ultra thick, red]
(9, 1.75) -- (9, 1) -- (10, 1) -- (10,0) -- (11,0) -- (11, -.25);

\draw[ultra thick, red]
(9.5, 1.75) -- (9.5, 1.5) -- (10.5, 1.5) -- (10.5, .5) -- (11.5, .5) -- (11.5, -.25);

\draw[ultra thick, red]
(11, 1.75) -- (11, 1) -- (11.5, 1) -- (11.5, 1.75);

\draw[ultra thick, blue, dashed]
(-.25, 1.5) -- (0, 1.5) -- (0, 1.75);
\draw[ultra thick, blue, dashed]
(-.25, 0) -- (0, 0) -- (0, -.25);
\draw[ultra thick, blue, dashed]
(-.25, 1) -- (1,1) -- (1, 0) -- (2, 0) -- (2, -.25);
\draw[ultra thick, blue, dashed]
(-.25, .5) -- (.5,.5) -- (.5, -.25);
\draw[ultra thick, blue, dashed]
(.5, 1.75) -- (.5, 1.5) -- (1.5, 1.5) -- (1.5, .5) -- (3, .5) -- (3, 1) -- (2, 1) -- (2, 1.75);
\draw[ultra thick, blue, dashed]
(2.5, -.25) -- (2.5, 0) -- (4, 0) -- (4, -.25);
\draw[ultra thick, blue, dashed]
(2.5, 1.75) -- (2.5, 1.5) -- (3.5, 1.5) -- (3.5, .5) -- (4.5, .5) -- (4.5, -.25);
\draw[ultra thick, blue, dashed]
(4, 1.75) -- (4, 1) -- (5, 1) -- (5, 0) -- (6, 0) -- (6, -.25);
\draw[ultra thick, blue, dashed]
(6, 1.75) -- (6, 1.5) -- (4.5, 1.5) -- (4.5, 1.75);

\draw[ultra thick, blue, dashed]
(5.5, .5) -- (5.5, 1) -- (7,1) -- (7, 0) -- (8, 0) -- (8, -.25);
\draw[ultra thick, blue, dashed]
(5.5, .5) -- (6.5,.5) -- (6.5, -.25);
\draw[ultra thick, blue, dashed]
(6.5, 1.75) -- (6.5, 1.5) -- (7.5, 1.5) -- (7.5, .5) -- (9, .5) -- (9, 1) -- (8, 1) -- (8, 1.75);
\draw[ultra thick, blue, dashed]
(8.5, -.25) -- (8.5, 0) -- (10, 0) -- (10, -.25);
\draw[ultra thick, blue, dashed]
(8.5, 1.75) -- (8.5, 1.5) -- (9.5, 1.5) -- (9.5, .5) -- (10.5, .5) -- (10.5, -.25);
\draw[ultra thick, blue, dashed]
(10, 1.75) -- (10, 1) -- (11, 1) -- (11, 0) -- (11.75, 0);
\draw[ultra thick, blue, dashed]
(11.75, 1.5) -- (10.5, 1.5) -- (10.5, 1.75);
\draw[ultra thick, blue, dashed]
(11.75, 1) -- (11.5, 1) -- (11.5, .5) -- (11.75, .5);

\foreach \x in {0, .5, 1, 1.5, 2, 2.5, 3, 3.5}{
    \foreach \y in {0,.5, 1,1.5}{
    \node[style={circle, draw, fill=black},scale=.4] at (\x,\y){};}};

\foreach \x in {4, 4.5, 5, 5.5, 6, 6.5, 7, 7.5}{
    \foreach \y in {0,.5,1,1.5}{
    \node[style={circle, draw, fill=black},scale=.4] at (\x,\y){};}};

\foreach \x in {8, 8.5, 9, 9.5, 10, 10.5, 11, 11.5}{
    \foreach \y in {0,.5,1,1.5}{
    \node[style={circle, draw, fill=black},scale=.4] at (\x,\y){};}};

\end{tikzpicture}

\vspace{.1in}

\begin{tikzpicture}[scale = 1.2]

\node at (-.5,1.5){0};
\node at (-.5,1){1};
\node at (-.5,.5){2};
\node at (-.5,0){3};

\node at (0,-.5){0};
\node at (1,-.5){2};
\node at (2,-.5){4};
\node at (3,-.5){6};
\node at (4,-.5){8};
\node at (5,-.5){10};
\node at (6,-.5){12};
\node at (7,-.5){14};
\node at (8,-.5){16};
\node at (9,-.5){18};
\node at (10,-.5){20};
\node at (11,-.5){22};

\draw[ultra thick, red] (1, -.25) -- (1,0) -- (0,0) -- (0,1.5) -- (.5,1.5) -- (.5,.5) -- (1.5, .5) -- (1.5, -.25);

\draw[ultra thick, red] (1, 1.75) -- (1, 1) -- (2, 1) -- (2,0) -- (3,0) -- (3,-.25);

\draw[ultra thick, red] (1.5, 1.75) -- (1.5, 1.5) -- (2.5, 1.5) -- (2.5, .5) -- (3.5, .5) -- (3.5, -.25);

\draw[ultra thick, red] (3, 1.75) -- (3, 1) -- (4, 1) -- (4,0) -- (5,0) -- (5,-.25);

\draw[ultra thick, red] (3.5, 1.75) -- (3.5, 1.5) -- (4.5, 1.5) -- (4.5, .5) -- (5.5, .5) -- (5.5, -.25);

\draw[ultra thick, red] (5, 1.75) -- (5, 1) -- (6, 1) -- (6,0) -- (7,0) -- (7,-.25);

\draw[ultra thick, red] (5.5, 1.75) -- (5.5, 1.5) -- (6.5, 1.5) -- (6.5, .5) -- (7.5, .5) -- (7.5, -.25);

\draw[ultra thick, red] (7, 1.75) -- (7, 1) -- (8, 1) -- (8,0) -- (9,0) -- (9,-.25);

\draw[ultra thick, red] (7.5, 1.75) -- (7.5, 1.5) -- (8.5, 1.5) -- (8.5, .5) -- (9.5, .5) -- (9.5, -.25);

\draw[ultra thick, red] (9, 1.75) -- (9, 1) -- (10, 1) -- (10,0) -- (11,0) -- (11,-.25);

\draw[ultra thick, red] (9.5, 1.75) -- (9.5, 1.5) -- (10.5, 1.5) -- (10.5, .5) -- (11.5, .5) -- (11.5, -.25);

\draw[ultra thick, red] (11, 1.75) -- (11, 1) -- (11.5, 1) -- (11.5, 1.75);

\draw[ultra thick, blue, dashed] (0,1.75) -- (0,1.5) -- (-.25, 1.5);

\draw[ultra thick, blue, dashed] (-.25,1) -- (1, 1) -- (1, 0) -- (2, 0) -- (2, -.25);

\draw[ultra thick, blue, dashed] (.5, 1.75) -- (.5, 1.5) -- (1.5, 1.5) -- (1.5, .5) -- (2.5, .5) -- (2.5, -.25);

\draw[ultra thick, blue, dashed] (2, 1.75) -- (2, 1) -- (3, 1) -- (3, 0) -- (4, 0) -- (4, -.25);

\draw[ultra thick, blue, dashed] (2.5, 1.75) -- (2.5, 1.5) -- (3.5, 1.5) -- (3.5, .5) -- (4.5, .5) -- (4.5, -.25);

\draw[ultra thick, blue, dashed] (4, 1.75) -- (4, 1) -- (5, 1) -- (5, 0) -- (6, 0) -- (6, -.25);

\draw[ultra thick, blue, dashed] (4.5, 1.75) -- (4.5, 1.5) -- (5.5, 1.5) -- (5.5, .5) -- (6.5, .5) -- (6.5, -.25);

\draw[ultra thick, blue, dashed] (6, 1.75) -- (6, 1) -- (7, 1) -- (7, 0) -- (8, 0) -- (8, -.25);

\draw[ultra thick, blue, dashed] (6.5, 1.75) -- (6.5, 1.5) -- (7.5, 1.5) -- (7.5, .5) -- (8.5, .5) -- (8.5, -.25);

\draw[ultra thick, blue, dashed] (8, 1.75) -- (8, 1) -- (9, 1) -- (9, 0) -- (10, 0) -- (10, -.25);

\draw[ultra thick, blue, dashed] (8.5, 1.75) -- (8.5, 1.5) -- (9.5, 1.5) -- (9.5, .5) -- (10.5, .5) -- (10.5, -.25);

\draw[ultra thick, blue, dashed] (10, 1.75) -- (10, 1) -- (11, 1) -- (11, 0) -- (11.75, 0);

\draw[ultra thick, blue, dashed] (10.5, 1.75) -- (10.5, 1.5) -- (11.75, 1.5);

\draw[ultra thick, blue, dashed] (11.75, 1) -- (11.5, 1) -- (11.5, .5) -- (11.75, .5);

\draw[ultra thick, blue, dashed] (-.25, .5) -- (.5, .5) -- (.5, -.25);

\draw[ultra thick, blue, dashed] (-.25, 0) -- (0,0) -- (0,-.25);

\foreach \x in {0, .5, 1, 1.5, 2, 2.5, 3, 3.5}{
    \foreach \y in {0,.5, 1,1.5}{
    \node[style={circle, draw, fill=black},scale=.4] at (\x,\y){};}};

\foreach \x in {4, 4.5, 5, 5.5, 6, 6.5, 7, 7.5}{
    \foreach \y in {0,.5,1,1.5}{
    \node[style={circle, draw, fill=black},scale=.4] at (\x,\y){};}};

\foreach \x in {8, 8.5, 9, 9.5, 10, 10.5, 11, 11.5}{
    \foreach \y in {0,.5,1,1.5}{
    \node[style={circle, draw, fill=black},scale=.4] at (\x,\y){};}};

\end{tikzpicture}
\captionof{figure}{Top: The partitionable decomposition of $C_{24} \sq C_4$ into cycles of length four from the proof of Lemma~\ref{8ell} ($\ell = 6$). Top middle: Potential recoloring locations $S_1, \ldots, S_{24}$. Bottom middle: A partitionable decomposition of $C_{24} \sq C_4$ into cycles of length 24. Since $n=6$, $S_i$ is recolored unless $i \in \{6, 12, 18, 24\}$.
Bottom: A partitionable decomposition of $C_{24} \sq C_4$ into cycles of length 96. Since $n=24$, $S_i$ is recolored unless $i =24$.
}
\label{C_24 sq C_4, recoloring & initial}

\end{figure}
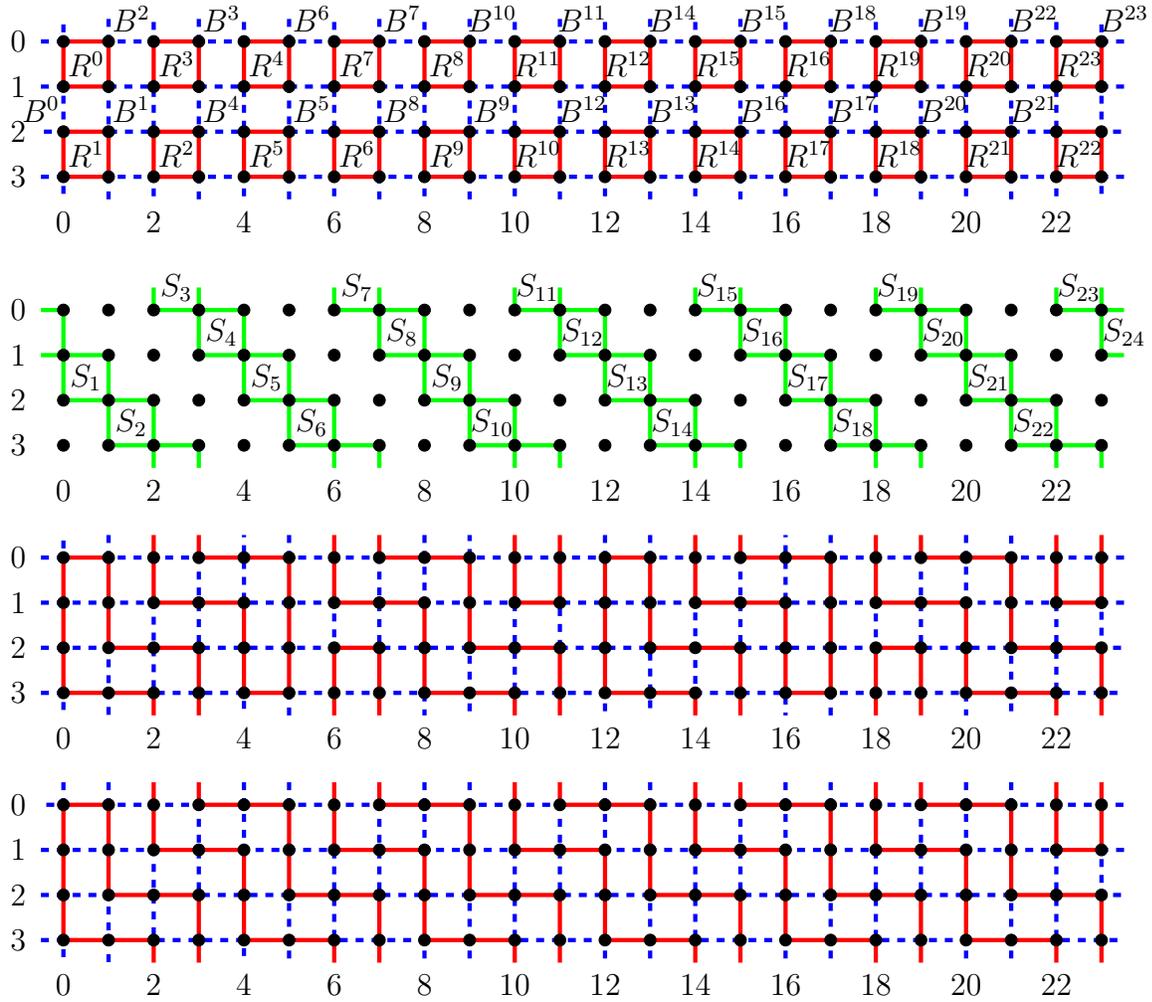

\begin{proposition}\label{cart}
For $\ell \ge 1$, if $G$ has a partitionable decomposition into cycles of length $4\ell$, then $G \sq C_4$ has a partitionable decomposition into cycles of length $4\ell$.
\end{proposition}

\begin{proof}
 We use the following notation. Let $V(C_4) = \{0, 1, 2, 3\}$ and  $E(C_4) = \{\{0,1\}, \{1,2\}, \{2,3\}, \{3,0\}\}$. Then
$V(G \sq C_4) = \{(v, i): v \in V(G), 0 \leq i \leq 3 \}$ and 
$E(G \sq C_4) = \{ (u, i)(v, j): u = v, \{i,j\} \in E(C_4) \text{ or } i = j, uv \in E(G)\}$. For any subgraph $H$ of $G$ and $0 \le j \le 3$, let $H^j$ be the subgraph of $G \sq C_4$ such that $V(H^j) = \{ (v, j) : v \in V(H)\}$ and $E(H^j) = \{(u,j)(v, j): uv \in E(H) \}$. Note that for all $j$, $H^j$ is isomorphic to $H$. Also note that the vertex sets of $G^j$ for $0 \leq j \leq 3$ partition the vertices of $G \sq C_4$. 

Suppose the decomposition of $G$ into cycles of length $4\ell$ has $k$ partition sets called $\F_1, \ldots, \F_k$. 
For $1 \le i \le k$, define the graph $F_i = \cup_{C \in \F_i} C$.  Since $\F_i$ is a partition set for $1 \le i\le k$, for all $i$, $V(F_i) = V(G)$, and $G = F_1 \sqcup \dots \sqcup F_k$ is a partitionable decomposition of $G$, where each partition set contains exactly one of the $F_i$.

For $1 \le i \le k-1$, let  $G_i = F_i^0 \cup F_i^1 \cup F_i^2 \cup F_i^3$, and let $G_k = F_k \sq C_4$. Note that since $V(F_k) = V(G)$, $G_k = F_k \sq C_4$ contains all vertical edges in $G \sq C_4$ (edges that differ in the second coordinate). Thus since $G = F_1 \sqcup \dots \sqcup F_k$, $G \sq C_4 = G_1 \sqcup \dots \sqcup G_{k-1} \sqcup G_k$, and for $1 \le i \le k$, $V(G_i) = V(G \sq C_4)$.  

Since for $1 \le i \le k$, $V(G_i) = V(G \sq C_4)$, by Proposition~\ref{ddd}, it remains to give a partitionable decomposition of $G_i$ into cycles of length $4\ell$ for $1 \le i \le k$.
For $1 \le i \le k-1$, and $0 \leq j \leq 3$, let $\F_{i,j} = \{C^j : C \in \F_i \}$, and let $\F'_i = \bigcup\limits_{0 \leq j \leq 3} \F_{i,j}$.  For $1 \le i \le k-1$, since the cycles in $\F_i$ form a partitionable decomposition of $F_i$ with one partition set, the cycles in $\F_i'$ form a decomposition of $G_i$ with one partition set.

We now describe a partitionable decomposition of $G_k = F_k \sq C_4$. 
Letting $G = F_k$, $\F = \F_k$, and $H=C_4$, Proposition~\ref{decomp} implies that $G_k$ can be decomposed into the union of tori $\sqcup_{C \in \F_k} C \sq C_4$.  Choosing $n = \ell$, Lemma~\ref{8ell} implies that each torus in the union has a partitionable decomposition into cycles of length $4n = 4\ell$ with two partition sets.  Suppose one partition set for each torus contains red cycles and the other blue.  Then letting $\F_k'$ be the union of the partition sets of red cycles and $\F_k''$ be the union of the partition sets of blue cycles, we have a partition of $G_k$ into cycles of length $4 \ell$ with two partition sets $\F_k'$, and $\F_k''$.

Putting it all together, we have a partitionable decomposition of $G \sq C_4$ into cycles of length $4\ell$ with partition sets $\F_1', \ldots \F_{k-1}', \F_k', \F_k''$.
\end{proof}

We are now ready to prove our main result.

\begin{theorem}\label{partitionable}
If $n$ is even, $n \ge 2$, and $2 \le i \le n$, then $Q_n$ has a partitionable decomposition into cycles of length $2^i$.
\end{theorem}

\begin{proof}

Fix $i\ge 2$.  We show by induction  on even $n$ that $Q_n$ has a partitionable decomposition into cycles of length $2^i$ for all even $n \ge i$.

Our base case depends on whether $i$ is even or odd.  If $i$ is odd, then the smallest even $n$ where $n \ge i$ is $n=i+1$.  In this case, Proposition~\ref{halfham} implies that $Q_n$ has a partitionable decomposition into cycles of length $2^{n-1} = 2^i$. If $i$ is even, then the smallest even $n$ where $n \ge i$ is  $n=i$.  In this case, Theorem~\ref{ham} implies that $Q_n$ has a partitionable decomposition into cycles of length $2^n = 2^i$. 

For the inductive step, assume for some even $n \ge 2$ that $Q_n$ has a partitionable decomposition into cycles of length $2^i$, and note that since $i \ge 2$, $2^i = 4\ell$ for some $\ell \ge 1$. Since $Q_{n+2} = Q_n \sq C_4$, letting $G=Q_n$,  Proposition~\ref{cart} guarantees that $Q_{n+2}$ has a partitionable decomposition into cycles of length $4\ell = 2^i$.
\end{proof}

\junk{
Induct on $n$.
BC: $n = 2$. Therefore, $2 \leq i \leq n \rightarrow i = 2$. Want a partitionable decomposition of $Q_2 = C_4$ into $C_4$, clearly works.
\\
\\
\textbf{IH: } Assume there is a partitionable decomposition for all $2 \leq i \leq n$ for $Q_n$.
\textbf{IS: } WTS there is a partitionable decomposition for all $2 \leq i \leq n + 2$. 
\\
Case 1: i = 2. Do separately.
\\
Case 2: $3 \leq i \leq n$. Then $2^i = 2^32^{i-3} = 8*2^{i-3}$, works by proposition 3.
\\
Case 3: $i = n+1$ or $i = n+2$, propositions 4,7.
}



\begin{thebibliography}{99}



\bibitem{alspach} B. Alspach, J. Bermond, and D. Sotteau. 
\newblock Decomposition into cycles I: Hamilton decompositions. In G. H. et al.\ (editors) Cycles and rays. Kluwer Academic (1990), 9--18.


\bibitem{AR}  D. Anick and M. Ramras. Edge decompositions of hypercubes by paths. Australas. J. Combin., 61 (2015), 210--226.

\bibitem{aubert}
J. Aubert  and B. Schneider. 
\newblock Decomposition de la somme cartesienne d'un cycle et de l'union de deux cycles hamiltoniens en cycles hamiltoniens. Discrete Mathematics 38.1 (1982), 7--16.




\bibitem{AOT21} M. Axenovich, D. Offner, and C. Tompkins. \newblock
Long path and cycle decompositions of even hypercubes. 
European J. Combin. 95 (2021), 103320, 20 pp.








\bibitem{erde} 
J. Erde. 
\newblock Decomposing the cube into paths. Discrete Math., 336 (2014), 41--45. 


\bibitem{fink}
J. Fink.
\newblock On the decomposition of $n$-cubes into isomorphic trees. J. Graph Theory 14 (1990), 405--411. 









\bibitem{HSW} P. Horak, J. Siran, and W. Wallis. \newblock Decomposing Cubes, J. Austral. Math. Soc. Ser. A 61.1 (1996), 119--128. 




\bibitem{kotzig}
A. Kotzig.
\newblock Every Cartesian product of two circuits is decomposable into two hamiltonian circuits,
Centre de Recherches Mathematique, Montreal, (1973).


\bibitem{MR}
M. Mollard and M. Ramras.
\newblock Edge decompositions of hypercubes by paths and by cycles. 
Graphs Combin. 31 (2015), no. 3, 729–-741.


\bibitem{ringel}
G. Ringel. 
\newblock \"Uber drei kombinatorische Probleme am $n$-dimensionalen W\"urfel und W\"urfelgitter.  Abhandlungen aus dem Mathematischen Seminar der Universit\"at Hamburg, Vol. 20, No. 1-2 (1955).



\bibitem{stout}
 Q. Stout.
 \newblock Packings in hypercubes, presented at 21st Southeastern Intl. Conf. on Combinatorics, Graph Theory, and Computing, Boca Raton, FL, (1990). Abstract available at \url{http://web.eecs.umich.edu/~qstout/abs/hyppack.html}.


\bibitem{TWB} S. Tapadia, B. Waphare, and  Y. Borse.  Cycle decompositions of the Cartesian product of cycles. 
Australas. J. Combin. 74 (2019), 443--459.

\bibitem{WW12}
S. Wagner and M. Wild.
\newblock Decomposing the hypercube $Q_n$ into $n$ isomorphic edge-disjoint trees. 
Discrete Math. 312 (2012), no. 10, 1819--1822.




\end{thebibliography}
\end{document}